\documentclass[12pt,draftcls,onecolumn]{IEEEtran}

\usepackage{amsmath,amssymb,amsfonts}

\usepackage[english]{babel}
\usepackage{algorithm,algorithmic,bm}
\usepackage{graphicx}
\usepackage{textcomp}
\usepackage{color}
\usepackage{lipsum}
\usepackage{epstopdf}
\usepackage{amsopn}
\usepackage{amsthm}
\usepackage{enumerate}
\usepackage{cancel}
\usepackage{mathrsfs}
\usepackage{mathdots}
\usepackage{euscript}
\usepackage{amscd}
\usepackage{cite}
\usepackage{placeins}
\usepackage[latin1]{inputenc}
\usepackage{tikz}
\usetikzlibrary{snakes,arrows,shapes}

\usepackage{subcaption}

\newtheorem{theorem}{Theorem}
\newtheorem{corollary}{Corollary}

\newtheorem{proposition}{Proposition}

\theoremstyle{definition}

\newtheorem{example}{Example}

\newtheorem{assumption}{Assumption}
\theoremstyle{remark}
\newcommand{\tr}{\mathrm{tr}}
\newcommand{\Gc}{\mathcal{G}}
\newcommand{\teo}{\begin{theorem}}
 \newcommand{\eteo}{\end{theorem}}
\newcommand{\prop}{\begin{proposition}}
 \newcommand{\eprop}{\end{proposition}}
\newcommand{\Ts}{\mathbb T}
\newcommand{\Bk}{\mathfrak B}
\newcommand{\Ss}{\mathbb{S}}
\newcommand{\al}[1]{\begin{align}#1\end{align}}
\newcommand{\su}[1]{  \begin{subequations} #1\end{subequations} }
\newcommand{\bmat}{\begin{bmatrix}}
\newcommand{\emat}{\end{bmatrix}}

\newcommand{\Rbb}{\mathbb R}

\newcommand{\Zbb}{\mathbb Z}
\newcommand{\Sbb}{\mathbb S}
\newcommand{\Nbb}{\mathbb N}
\newcommand{\Tbb}{\mathbb T}

\newcommand{\pb}{\mathbf  p}

\newcommand{\Eb}{\mathbf E}

\newcommand{\Gb}{\mathbf G}

\newcommand{\Qb}{\mathbf Q}

\newcommand{\Wb}{\mathbf W}
\newcommand{\Sc}{\mathcal S}

\newcommand{\Ls}{\mathbb{L}}

\newcommand{\nn}{\nonumber}

\begin{document}

\title{On the identification of ARMA graphical models}

\author{Mattia Zorzi
\thanks{M.~Zorzi is with the Department of Information Engineering, University of Padova, Via Giovanni Gradenigo, 6b, 35131 Padova, Italy (email: \texttt{zorzimat@dei.unipd.it}).}
}

\maketitle

\begin{abstract}
The paper considers the problem to estimate a graphical model corresponding to an autoregressive moving-average (ARMA) Gaussian stochastic process. We propose a new maximum entropy covariance and cepstral extension problem and we show that the problem admits an approximate solution which represents an ARMA graphical model whose topology is determined by the selected entries of the covariance lags considered in the extension problem. Then, we show how the corresponding dual problem is connected with the maximum likelihood principle. Such connection allows to design a Bayesian model and characterize an approximate maximum a posteriori estimator of the ARMA graphical model in the case the graph topology is unknown. We test the performance of the proposed method through some numerical experiments. 
\end{abstract}


\section{Introduction}\label{sec:intro}
 
In various scenarios, we gather data of high dimensionality that characterizes intricate systems with a multitude of variables. Consequently, there is the necessity to comprehend the interrelationships among these variables. Essentially, we are compelled to extract from the data the underlying structure of the graphical model that depicts these associations. Clearly, a graphical model holds significance when it exhibits a certain degree of sparsity, meaning it is not densely connected in terms of edges. Moreover, the network structure depends on how we define these interactions. The latter can describe conditional dependence relations between the stochastic processes characterizing the variables of the system and whose resulting network is undirected 
 \cite{e20010076,you2023sparse,LATENTG,KRON}. Alternatively, the interactions can describe whether a variable of the system is necessary to reconstruct another variable and in this case the resulting network is directed   \cite{MATERASSI2,9600870,SHAIKHVEEDU2023111182,doddi2020estimating}. It is noteworthy that a related topic regards the identification of network modules solely from neighboring node data, as discussed in works like \cite{HOF1,HOF2}.

In this paper we focus our attention on the identification of graphical models whose edges encode conditional dependence relations. The typical approach is to set up a multivariate regularized spectral estimation problem which originates from a maximum entropy covariance extension problem. Many generalizations of the covariance extension paradigm have been proposed in the classic scenario, i.e. the case in which the spectrum is not associated to a graphical model, \cite{FPR-08,RFP-09,FRT-11,zhu2018wellposed,zhu2020m}. Interestingly, in the case one want to infer the moving-average (MA) part of the model, then it is sufficient to consider a maximum entropy covariance extension problem to which we add some moments constraints regarding the cepstral coefficients of the spectrum, see for instance the works  \cite{ringh2018multidimensional,KLR-16multidimensional,RKL-16multidimensional,enqvist2004aconvex}.

Most of the identification paradigms for graphical models restrict the attention to autoregressive (AR) processes, see for instance \cite{SONGSIRI_TOP_SEL_2010,MAANAN2017122,ALPAGO_SL_REC}. The reason is that the parametrization of the spectrum induced by this model class allows to impose easily sparsity constraints in terms of edges of the graphical model.

Learning ARMA graphical models has been a less-explored problem in the literature. The first paper, \cite{6365751}, fixed the MA part using a scalar trigonometric  polynomial, which is estimated through a preliminary step, then a maximum entropy covariance extension problem is considered. A generalization has been proposed in \cite{LINKPRED} where the MA part is fixed through a matrix trigonometric polynomial obeying a certain structure. Recently, an iterative procedure to estimate the AR and MA parts has been proposed in \cite{YOU2022110319}. More precisely, the estimation of the AR and MA parts is preformed separately according two different paradigms. In addition, the sparsity of the graph is induced through a group Lasso regularization  term which can introduce an adverse bias in general, see \cite{WIPF_2010}.

In this paper we propose a new maximum entropy covariance and cepstral extension problem, i.e. the optimal spectrum maximizes the entropy rate while matching some cepstral coefficients and some entries of some 
covariance lags. We show that the problem admits an approximate solution which represents an ARMA graphical model whose topology is determined by the selected entries of the covariance lags considered in the extension problem. The fact that the existence result regards only the approximate solution is in line with the results found in univariate case, see for instance \cite{RKL-16multidimensional}. It is worth noting that in this paradigm the AR and MA parts are jointly estimated. Then, we show how the corresponding dual problem is connected with the maximum likelihood principle. The latter allows to design a Bayesian model and characterize an approximate maximum a posteriori estimator of the ARMA graphical model in the case the graph topology is unknown. In particular, 
the sparsity of the graph is induced through a reweighted iterative scheme, earlier used for the AR case in \cite{ZORZI2019108516}, in which the weights are iteratively adapted in such a way to reduce the presence of the adverse bias.

The outline of the paper is as follows. In Section \ref{sec:back} we provide the necessary background about dynamic graphical models. In Section \ref{sec:extension} we introduce the maximum entropy covariance and cepstral extension problem, while in Section \ref{sec_dualpb} we show its regularized dual problem admits solution. In Section \ref{sec:ML} we show how the dual problem is connected with the maximum likelihood estimator. In Section \ref{sec:Bayes} we propose a Bayesian model to estimate an ARMA graphical model, including its topology, from the data. Section \ref{sec:example} contains some numerical examples showing the performance of the proposed method. Finally, we draw the conclusions in Section \ref{sec:concl}.

{\em Notation.} $\Nbb$ and $\Zbb$ are the set of natural and integer numbers, respectively; $\Rbb$ and $\Rbb_+$ is the set of real and positive real numbers, respectively. $\Sc^m$ denotes the vector space of symmetric matrices of dimension $m$, $\Sc_+^m$ is the cone of positive definite symmetric matrices of dimension $m$. $\Ls_\infty^m(\Ts)$ denotes  the space of functions defined on $\Ts:=[0,2\pi)$, taking values in the space of Hermitian matrices of dimension $m$ and bounded  a.e.  Let $\Ss^m_+\subset \Ls^m_\infty(\Ts)$ denote the family of spectral densities on $\Ts$ of dimension $m$ which are bounded and coercive a.e.. Given $\Phi(e^{i\theta})\in\Ls^m_\infty(\Ts)$, $\int_\Ts\Phi$ denotes the integral of $\Phi$ over $\Ts$ with respect to the normalized Lebesgue measure $\mathrm d \theta/2\pi$. Given $X\in\Rbb^{m\times m}$, $[X]_{jh}$ denotes the entry of $X$ in position $(j,h)$; this notation applies also for $\Phi\in\Ls^m_\infty(\Ts)$. We denote as $\mathrm{supp}(\Phi)$ the support of $\Phi$, i.e. it is a binary matrix whose entries are equal to one if the corresponding entries in $\Phi$ are different from the null function.

\section{Background on graphical models}\label{sec:back}
A dynamic undirected graphical model is a graph ${\cal G} (V,E)$ associated to a  $m$-dimensional Gaussian, stationary, purely non-deterministic stochastic process $y:=\{y(t),\; t\in\Zbb\}$ having vertexes $V:=\{1\ldots m\}$ corresponding to the components of $y$ and edges $E\subseteq V\times V$ describing the conditional dependence relations between the components. More precisely, for any index set $I\subseteq V$, define ${\cal X}_I:=\overline{\mathrm{span}}\{ y_j(t) \hbox{ s.t. } j\in I, \; t\in\Zbb\}$ the closure of the set containing all the finite linear combinations of the random variables $y_j(t)$, $j\in I$, so that for any $j\neq h$, the notation 
\al{ {\cal X}_j\; \bot\; {\cal X}_h\, |\,{\cal X}_{V\setminus \{j,h\}}\nn} 
means that for all $t,s$, $y_j(t)$ and $y_h(s)$ are conditionally independent given the space linearly generated by $\{ y_k(t), \; k\in V\setminus \{j,h\}, \; t\in\Zbb\}$. Let $\Phi(e^{i\theta})$, $\theta\in\Tbb$, denote the power spectral density of $y$. Then, it is  possible to prove that, \cite{brillinger1996}:
\al{{\cal X}_j\; \bot\; {\cal X}_h\, |\,{\cal X}_{V\setminus \{j,h\}}\quad \iff \quad [\Phi(e^{i\theta})^{-1}]_{jh}=0\nn}
and thus
\al{(j,h)\notin E\quad \iff \quad [\Phi(e^{i\theta})^{-1}]_{jh}\equiv 0.\nn}
In plain words, the topology of the graph corresponds to the support of the inverse of the spectral density. Notice that, $(j,j)\in E$ for any $j\in V$.

\section{Covariance and Cepstral extension for graphical models}\label{sec:extension}
Consider a $m$-dimensional Gaussian stationary stochastic process $y$ with power spectral density $\Phi$. Let $R_k:=\mathbb E [y_{t+k}y_{t}^T]$ be the $k$-th covariance lag of $y$. Moreover, we define as
\al{\label{def_cepstral}c_k:=\int_{\Tbb}e^{i\theta k}\log \det\Phi}
 the $k$-th cepstral coefficient of $y$, \cite{lindquist2013multivariate}.  
Assume to know $[R_k]_{jh},c_k\in \mathbb R$ with $k=0\ldots n$ and $(j,h)\in E$. Then, we consider the following maximum entropy covariance and cepstral extension problem:
\su{\label{moment_pb} \al{& \underset{\Phi\in\Ss_+^m}{ \max\;} \int_\Ts \log\det \Phi\\ & \label{cov_c} \hbox{ s.t. }  \int_\Ts e^{i\theta k}  [\Phi]_{jh}=[R_k]_{jh},\quad k=0 \ldots n,\; (j,h)\in E \\
 & \label{cepstral_c} \hspace{0.5cm} \int_\Ts e^{i\theta k}   \log \det\Phi=c_k,\quad 1 \ldots n.
} } 
It is worth noting that (\ref{moment_pb}) is a generalization of the covariance selection problem proposed by Dempster, \cite{DEMPSTER_1972}.  
The Lagrangian is 
 \al{L(\Phi,& \pb,\Qb)=  \int_\Ts\log\det \Phi\nn\\ &-   \sum_{(j,h)\in E}\sum_{k=0}^n [Q_k]_{jh}   \left(\int_\Ts e^{i\theta k}  [\Phi]_{jh} -[R_k]_{jh}\right) \nn\\
 &+\sum_{k=1}^n p_k\left(\int_\Ts e^{i\theta k}   \log\det \Phi-c_k\right)-p_0c_0\nn \\
=&   \int_\Ts \pb \log \det\Phi-\tr(\Qb\Phi) +\sum_{k=0}^n\tr( Q_k^T R_k)-\sum_{\substack{k=0}}^n p_k c_k \nn  
}
 where $p_0=1$, $c_0$ is the zeroth cepstral coefficient  defined as in (\ref{def_cepstral}) with $k=0$ and the constant term $p_0c_0$ is only added for convenience in the following;  $Q_k$, with $[Q_k]_{jh}=0$ for $(j,h)\notin E$, and $p_k$, with $k\neq 0$, are the Lagrange multipliers and 
 \al{ \label{defpQ}\pb(e^{i\theta})&:=1+\frac{1}{2}\sum_{\substack{k=1}}^n p_k (e^{-i\theta k}+ e^{i\theta k}),\; \nn\\
 \Qb(e^{i\theta})&:=Q_0+\frac{1}{2}\sum_{k=1}^n Q_k e^{-i\theta k}+Q_k^T e^{i\theta k}
  }
 are the corresponding trigonometric polynomials. Let 
 {\scriptsize \al{\Bk_+^m(E)&:=\left\{\, Q(e^{i\theta})=Q_0+\frac{1}{2}\sum_{k=1}^n Q_k e^{-i\theta k}+Q_k^T e^{i\theta k} \hbox{ s.t.} \right. \nn \\ & \hspace{1.8cm}\left. \,[Q_k]_{jh}=0 \hbox{ for } (j,h)\notin E, Q(e^{i\theta})>0 \hbox{ on } \Ts\right\}\nn\\
 \Bk_+^o&:=\left\{\, p(e^{i\theta})=1+\frac{1}{2}\sum_{k=1}^n p_k (e^{-i\theta k}+e^{i\theta k})\hbox{ s.t. }p(e^{i\theta})>0 \hbox{ on } \Ts\right\}.\nn }}Let $\overline{\Bk_+^o}$ denote the closure of ${\Bk_+^o}$; $\overline{\Bk_+^m}(E)$ denotes the closure of $\Bk_+^m(E)$ without the points $\Qb$ for which $\det(\Qb(e^{i\theta}))=0$ for all $\theta\in\Tbb$. The dual function $\sup_{\Phi \in \Sbb_+^m}\, L(\Phi,\pb,\Qb)$ is finite if and only if $\pb \in\overline{ \Bk_+^o}$ and $\Qb\in\overline{ \Bk_+^m}(E)$. Indeed,  assume that there exists $\theta_0$ for which $\pb(e^{i\theta_0})<0$. Taking $\Phi$ which tends to be singular in a neighborhood of $\theta_0$ we have that $L(\Phi,\pb,\Qb)\rightarrow \infty$. Hence, we must have $\pb\in\overline{\Bk^o_+}$.  Moreover, assume that there exists $\theta_0$ for which $\Qb(e^{i\theta_0})$ has at least one negative eigenvalue. Let $v$ be the corresponding eigenvector. Taking $\Phi$ which diverges along $v$ in a neighborhood of $\theta_0$ we have that $L(\Phi,\pb,\Qb)\rightarrow \infty$. Finally, if $\det(\Qb(e^{i\theta}))=0$ for all $\theta\in\Tbb$ then we can exploit the fact that there exists $\theta_0$ such that $\pb(e^{i\theta})>0$ in a neighborhood of $\theta_0$ because $\pb\in\overline{\Bk_+^o}$. Let $v(e^{i\theta})$ be a nonnull vector belonging to the null space of $\Qb(e^{i\theta})$. Taking $\Phi$ which diverges along $v(e^{i\theta})$ in the  neighborhood of $\theta_0$  we have that $L(\Phi,\pb,\Qb)\rightarrow \infty$. Thus, we must have $\Qb\in\overline{\Bk^m_+}(E)$. 
 
 Let $\pb\in \overline{\Bk_+^o}$ and $\Qb\in \overline{\Bk_+^m}(E)$, then $L(\cdot, \pb,\Qb)$ is concave on $\Ss_+^m$ and  the unconstrained maximum is also a stationary point, i.e. it is obtained by setting equal to zero the first variation of $L$ along the direction $\delta \Phi\in \Ls_\infty^m(\Ts)$:
 \al{\delta L &(\Phi, \pb, \Qb;\delta\Phi)= \tr  \int_\Ts  \pb \Phi^{-1}\delta\Phi  -\Qb \delta \Phi.\nn
 } Accordingly, the stationary point satisfies the condition: 
 \al{\label{opt_form}\Phi^\circ(e^{i\theta})=\pb(e^{i\theta}) \Qb(e^{i\theta})^{-1 } \hbox{ a.e. }. } 
The dual function takes the form  
{\small \al{J&(\pb,\Qb):=L(\Phi^\circ,\pb,\Qb)	\nn\\
&=\int_\Ts  \pb\log\det(\pb\Qb^{-1})-\pb m+\sum_{k=0}^n\tr( Q_k^T R_k)-\sum_{\substack{k=0}}^n p_k c_k.\nn
}}
 Therefore, the dual problem is 
 \al{\label{dual_pb}& \underset{\pb,\Qb}{ \min\;} J(\pb,\Qb)\\ & \hbox{ s.t. }  \pb\in \overline{\Bk_+^o},\; \Qb\in\overline{\Bk_+^m}(E). \nn   }
 If the dual problem admits an interior point solution, then Problem (\ref{moment_pb}) admits a unique solution of the form in (\ref{opt_form}). The latter represents the power spectral density of an ARMA process whose support of the inverse of its spectral density corresponds to $E$. Accordingly, model (\ref{opt_form}) describes an ARMA  graphical model whose conditional dependence relations are specified by $E$.
 At this point it is worth comparing the extension problem (\ref{moment_pb}) with the one proposed in \cite{6365751}; more precisely, the following maximum entropy covariance extension problem has been proposed
\al{& \underset{\Phi\in\Ss_+^m}{ \max\;} \int_\Ts \pb \log\det \Phi\nn\\ & \label{pb_avventi} \hbox{ s.t. }  \int_\Ts e^{i\theta k}  [\Phi]_{jh}=[R_k]_{jh},\quad k=0 \ldots n,\; (j,h)\in E
}
where $\pb(e^{i\theta })=\sum_{k=-n}^n p_k e^{-i\theta k}>0$ is given, e.g. it is usually estimated from the data though a preliminary step. It turns out that the solution to (\ref{pb_avventi})
has the same form of (\ref{opt_form}) and thus it corresponds to an ARMA graphical model, but the main difference is that $\pb$ is fixed by the user. On the contrary, our formulation is more general in the sense that $\pb$ and $\Qb$ are jointly optimized. 

  However, to prove the existence 
  of  an interior point solution to (\ref{dual_pb}) is not trivial, see the scalar case (i.e. $m=1$) in \cite{KLR-16multidimensional,Zhu-Zorzi2023cepstral}. Drawing inspiration from the scalar case \cite{enqvist2004aconvex}, in the next section we will address this issue considering a regularized version of the dual problem and whose regularizer forces the optimal solution $(\pb,\Qb)\in\Bk^o_+\times \Bk^m_+(E)$.

In the case one want to consider the general case in which the trigonometric polynomials $\pb$ and $\Qb$ have different degree, say $n_p$ and $n_q$, then it is sufficient to replace $n$ with $n_p$ and $n_q$ in (\ref{cepstral_c}) and (\ref{cov_c}), respectively. For ease of exposition we consider the case $n=n_p=n_q$, but all the following results trivially hold also for the case $n_p\neq n_q$.  
 
\section{Existence of the approximate solution} \label{sec_dualpb}
We consider the following regularized dual problem
 \al{\label{dual_pb_r}& \underset{\pb,\Qb}{ \min\;} J(\pb,\Qb)+g_\lambda(\pb) \\ & \hbox{ s.t. }  \pb\in \overline{\Bk_+^o},\; \Qb\in\overline{\Bk_+^m}(E), \nn }
 where \al{g_\lambda(\pb):=\lambda \int_\Ts \frac{1}{\pb}\nn}
is the regularization term and  $\lambda>0$ is the regularization parameter which is taken sufficiently small.
 The aim of this section is to prove the following result showing that Problem (\ref{dual_pb_r}) admits a unique solution belonging to the interior of the feasible set. 
 \begin{assumption}\label{hpPhi}There exists a power spectral density $\bar\Phi$ which is positive definite on some open interval in $\Tbb$ and  such that
 \al{\int_\Ts [\bar\Phi]_{jh}e^{i\theta k}=[R_k]_{jh}, \quad (j,h)\in E, \;k=0\ldots n. \nn}
\end{assumption} 
 \teo \label{teo_Demp_gen}  Under Assumption \ref{hpPhi}, Problem (\ref{dual_pb_r}) admits a unique solution $(\hat \pb,\hat \Qb)\in \Bk_+^o\times \Bk_+^m(E)$ and the rational spectral density $\hat \Phi =\hat \pb\hat \Qb^{-1}\in \Ss_+^m$ is such that 
 \su{\label{optcond}\al{ &\int_\Ts e^{i\theta k}  [\Phi]_{jh}=[R_k]_{jh},\quad k=0 \ldots n,\; (j,h)\in E\nn \\
 &  \int_\Ts e^{i\theta k}   \log \det\Phi=c_k+\varepsilon_k,\quad 1 \ldots n\nn}}
 where $\varepsilon_k=\lambda\int_\Ts e^{i\theta k}\pb^{-2}$, i.e. the cepstral coefficients are matched only approximately and the approximation error is 
$\varepsilon_k$ which depends on $\lambda$. \eteo 
It is worth noting that  $\hat \Phi$ represents an approximate solution of Problem (\ref{moment_pb}) and it corresponds to an ARMA graphical model whose graph is $\Gc(V,E)$. Let \al{\label{dualfreg} \tilde J_\lambda(\pb,\Qb):= J(\pb,\Qb)+g_\lambda(\pb)}
be the regularized dual function in Problem (\ref{dual_pb_r}). 
\prop \label{pro1}$\tilde J_\lambda$ is strictly convex on $\Bk^o_+\times \Bk^m_+(E)$.\eprop 
\IEEEproof Let $\Sc^m(E):=\{\,Y\in\Sc^m \hbox{ s.t. } [Y]_{jh}=0 \hbox{ for } (j,h)\notin E\}$ and $\Sc^m_+(E):=\Sc^m(E) \cap \Sc^m_+$. Consider the function
\al{h_\lambda(X,Y):= x\log\det(xY^{-1})+\lambda x^{-1} \nn}
with $x\in\Rbb_+$ and $Y\in \Sc_+^m(E)$. Then, its first and second variations along the direction $(\delta x, \delta Y)\in\Rbb\times \Sc^m(E)$ are  
\al{\delta & h_\lambda(x,Y; \delta x, \delta Y)= \nn\\ &(m+ m\log x-\log\det(Y)-\lambda x^{-2})\delta x- x\tr(Y^{-1}\delta Y) \nn}
\al{\delta^2  & h_\lambda(x,Y; \delta x, \delta Y)=(m x^{-1}+2 \lambda x^{-3})\delta x^2 \nn\\ &\hspace{1cm}- 2\tr(Y^{-1}\delta Y)\delta x+x \tr(Y^{-1}\delta YY^{-1}\delta Y) \nn \\&=\left[\begin{array}{cc}\delta x &  \mathrm{vec}(\delta Y)^T \end{array}\right]M
\left[\begin{array}{c}\delta x  \\\mathrm{vec}(\delta Y)  \end{array}\right]\nn
}
where \al{M:=\left[\begin{array}{cc} m x^{-1}+2 \lambda x^{-3} &\mathrm{vec}(Y^{-1})^T  \\\mathrm{vec}(Y^{-1})  & x(Y^{-1}\otimes Y^{-1}) \end{array}\right].\nn}
Matrix $M$ is positive definite because $m x^{-1}+2 \lambda x^{-3}>0$ and the Schur complement $M/ m x^{-1}+2 \lambda x^{-3} $ is positive:
\al{m &x^{-1}+2 \lambda x^{-3} -x^{-1} \mathrm{vec}(Y^{-1})^T(Y^{-1}\otimes Y^{-1})^{-1}\mathrm{vec}(Y^{-1})\nn\\
& = mx^{-1}+2 \lambda x^{-3} -x^{-1}\mathrm{vec}(Y^{-1})^T(Y\otimes Y)\mathrm{vec}(Y^{-1})\nn\\
& = mx^{-1}+2 \lambda x^{-3} -x^{-1}\mathrm{vec}(Y^{-1})^T\mathrm{vec}(Y)\nn\\
& = mx^{-1}+2 \lambda x^{-3} -x^{-1}\tr(Y^{-1}Y)=2 \lambda x^{-3} >0.\nn
}
We conclude that $h_\lambda(x,Y)$ is strictly convex in the domain $\Rbb_+\times \Sc^m_+(E)$ because $\delta^2  h_\lambda(x,Y; \delta x, \delta Y)>0$ for any nonnull direction $(\delta x, \delta Y)\in \Rbb\times \Sc^m(E)$. Finally, notice that 
\al{\tilde J_\lambda(\pb,\Qb):= \int_\Tbb h_\lambda(\pb,\Qb)-\pb m+\sum_{k=0}^n\tr( Q_k^T R_k)-\sum_{\substack{k=0}}^n p_k c_k\nn}
where the first term is strictly convex on $\Bk^o_+\times \Bk_+^m(E)$ because it is the integral of a strictly convex function while the remaining terms are linear in $(\pb,\Qb)$. Accordingly, the claim holds. 
\qed\\

In what follows we denote by $\partial \Bk_+^m(E)$ the boundary of $ \Bk_+^m(E)$ which is defined without the points $\Qb$ for which $\det(\Qb(e^{i\theta}))=0$ for all $\theta\in\Tbb$.
Moreover, we denote by $\partial \Bk_+^o$ the boundary of $\Bk_+^o$.
 
\prop \label{prop2}$\tilde J_\lambda$ is lower semicontinuous on $\overline{\Bk^o_+}\times \overline{\Bk^m_+}(E)$ with values on the extended reals. \eprop
\IEEEproof From Proposition \ref{pro1} it follows that $\tilde J_\lambda$ is differentiable and thus continuous in the interior of $\overline{\Bk^o_+}\times \overline{\Bk^m_+}(E)$. Accordingly, we only need to show the lower semicontinuity on the boundary of $\overline{\Bk^o_+}\times \overline{\Bk^m_+}(E)$. More precisely, we only need to show that the three terms 
\al{\tilde J_1=-\int_\Tbb \pb\log \det(\Qb),\;  \tilde J_2=\int_\Tbb m\pb \log\pb,\; \tilde J_3=\int_\Tbb \lambda \pb^{-1}\nn}
are lower semicontinuous because the other terms in $\tilde J_\lambda$ are obviously continuous. \\
{\em First term.} Take a sequence $\{(\pb_j,\Qb_j)\in\overline{\Bk^o_+}\times \overline{\Bk^m_+}(E), \, j\in\Nbb\}$ which converges to $(\pb,\Qb)$ and the latter belongs to the boundary of the feasible set.  Assume that $\Qb\in\partial \Bk^m_+(E)$ and $\pb\in\overline{\Bk^o_+}$. Since $\Qb$ is bounded and the convergence $\Qb_j\rightarrow \Qb$ is uniform, then we have that $\mu:=\sup_j \max_\theta \sigma_{max} (\Qb_j)<\infty$ where $\sigma_{max}(Y)$ denotes the maximum eigenvalue of $Y\in\Sc_+^m$. Accordingly, $0\leq \mu^{-1} \Qb_j \leq I$ for any $j\in\Nbb$. Let 
$\mathcal Z(\Qb):=\{\, \theta\in\Tbb \hbox{ s.t. } \det(\Qb(e^{i\theta}))=0\,\}$ and $\tilde{\mathcal Z}:=(\cup_{j\in\Nbb} \mathcal Z(\Qb_j) )\cup \mathcal Z(\Qb)$. Since $\tilde{\mathcal Z}$ has Lebesgue measure equal to zero, we have  \al{\tilde J_1= \int_\Tbb  \pb_j\log \det(\Qb_j) = 
\int_{\Tbb\setminus \tilde{\mathcal Z}}  \pb_j\log \det(\Qb_j),\nn}
thus we can consider the integrands over $\Tbb\setminus \tilde{\mathcal Z}$ in which they are well defined. Then, 
\al{\pb_j \log\det(\mu^{-1}\Qb_j)\rightarrow \pb\log\det(\mu^{-1}\Qb)\nn}
pointwise in $\Tbb\setminus \tilde{\mathcal Z}$ as $j\rightarrow \infty$. 
Since $-\pb_j \log\det(\mu^{-1}\Qb_j)$ is a nonnegative function, by Fatou's lemma \cite[p. 23]{rudin1987real} we have that
\al{\int_{\Tbb\setminus \tilde{\mathcal Z}}\pb \log\det(\mu^{-1}\Qb) \leq \underset{j\rightarrow \infty}{\mathrm{lim\, inf}}\int_{\Tbb\setminus \tilde{\mathcal Z}} \pb_j \log\det(\mu^{-1}\Qb_j)\nn}
which proves that $\tilde J_1$ is lower semicontinuous in these points. In the remaining case $\Qb\in\Bk^m_+(E)$,  $\pb\in \partial\Bk^o_+$ we have that $\tilde J_1$ is continuous in these points.\\
{\em Other terms.} The lower semicontinuity of $\tilde J_2$ has been proved in the proof of Lemma 5.1 in \cite{RKL-16multidimensional}. The lower semicontinuity of $\tilde J_3$ follows from the proof of Lemma 5.5 (first case) in \cite{Zhu-Zorzi2023cepstral}. \qed\\
 
The aim of the next propositions is to analyze the behavior of $\tilde J_\lambda$ when the optimization variables diverge or are on the boundary of the feasible set. We define the norm of the Lagrange multipliers:
\al{\|\pb\|:= \sqrt{\sum_{k=1}^n  p_k^2}, \quad \|\Qb\|:=\sqrt{\sum_{k=0}^n  \tr(Q_kQ_k^T)} \nn}
and $\|(\pb,\Qb)\|= \|\pb\|+\|\Qb\|$. 
 
 \prop \label{prop3} Let Assumption \ref{hpPhi} hold. Consider a sequence $\{(\pb_j,\Qb_j)\in\overline{\Bk^o_+}\times \overline{\Bk^m_+}(E), \, j\in\Nbb\}$ such that $\|(\pb_j,\Qb_j)\|\rightarrow \infty$ as $j\rightarrow\infty$. Then $\tilde J_\lambda(\pb_j,\Qb_j)\rightarrow \infty$.\eprop
 \IEEEproof Since $\pb_j$ is a nonnegative function and its integral on $\Tbb$ is equal to one, then it follows that $\|\pb_j\|$ cannot diverge as $j\rightarrow \infty$. Accordingly, it must happen that $\|\Qb_j\|\rightarrow\infty$. By Assumption \ref{hpPhi}, we have that
 \al{\tilde J_\lambda(\pb,\Qb) &\geq \int_\Tbb -\pb\log \det(\Qb)+\sum_{k=0}^n\tr(Q_kR_k^T)+s(\pb)\nn\\
 &\label{J_perse}=\int_\Tbb -\pb\log \det(\Qb)+\tr (\bar\Phi\Qb)+s(\pb)}
 where \al{s(\pb):=\int_\Tbb m(\pb\log \pb -\pb )-\sum_{k=0}^n p_kc_k. \nn} Accordingly, $s(\pb_j)$ cannot diverge as $j\rightarrow \infty$. Let $\Qb_j^0:=\Qb_j/\|\Qb_j\|\in\overline{\Bk_+^m}(E)$. Let \al{\eta:=\underset{j\rightarrow \infty}{\mathrm{\lim\,\inf }}\int_\Tbb \tr(\bar\Phi \Qb_j)\geq 0\nn .} Then, the sequence $\{ \Qb_{j},\; j\in\Nbb\}$ has a subsequence $\{ \Qb_{j_m},\; j_m\in\Nbb\}$ such that $\int_\Tbb \tr(\bar\Phi \Qb_{j_m})\rightarrow \eta$ as $j_m\rightarrow \infty$. Moreover, since such subsequence belongs to the surface of the unit ball, which is a compact set, it has a subsubsequence $\{ \Qb_{j_t},\; j_t\in\Nbb\}$ which is convergent. Let 
 \al{\Qb_\infty^0=\lim_{j_t\rightarrow \infty} \Qb_{j_t}\in\overline{\Bk_+^m}(E).\nn} Then, $\eta=\int_\Tbb \tr(\bar\Phi\Qb_\infty^0)$ and $\eta>0$. The latter property can be proved by contradiction. Assume that $\eta=0$. Since $\tr(\bar\Phi\Qb_\infty^0)\geq 0$ a.e. on $\Tbb$, it follows that $\tr(\bar\Phi\Qb_\infty^0)= 0$ a.e. on $\Tbb$. By Assumption \ref{hpPhi}, then $\Qb_\infty^0$  must vanish on some open interval in $\Tbb$. By \cite[Lemma 1]{ringh2015multidimensional}, it follows that $\Qb_\infty^0\equiv 0$ which is a contradiction because $\|\Qb_\infty^0\|=1$. Taking into account (\ref{J_perse}), we have
 \al{&\underset{j\rightarrow \infty}{\mathrm{lim\, \inf}} \tilde J_\lambda(\pb_j,\Qb_j)\geq \underset{j\rightarrow \infty}{\mathrm{lim\, \inf}} \,[\,\|\Qb_j\| \int_\Tbb \tr (\bar\Phi\Qb_j^0)\nn\\ &\hspace{0.8cm} -m\log\|\Qb_j\|-\int_\Tbb\pb_j\log\det(\Qb_j^0)+s(\pb_j)]\nn\\
 &=\underset{j\rightarrow \infty}{\mathrm{lim\, \inf}}[\,\eta\|\Qb_j\| -m\log\|\Qb_j\|]\nn\\ &\hspace{0.8cm}+\underset{j\rightarrow \infty}{\mathrm{lim\, \inf}}\,[ s(\pb_j) -\int_\Tbb\pb_j\log\det(\Qb_j^0)]=\infty\nn}
 where we exploited the fact that the second term does not diverge as $j\rightarrow \infty$.
 \qed\\
 
 \prop \label{prop4} Consider a sequence $\{(\pb_j,\Qb_j)\in \Bk^o_+ \times  \Bk^m_+ (E), \, j\in\Nbb\}$ converging to $(\pb,\Qb)$ which belongs to the boundary of $\overline{\Bk^o_+}\times \overline{\Bk^m_+}(E)$. Then, it cannot be an infimizing sequence for $\tilde J_\lambda$.\eprop
 \IEEEproof
 We have two cases. \\
 {\em First case: $\pb\in \partial \Bk^o_+$ and $\Qb\in \overline{\Bk^m_+}(E)$.} Let $Q_{k,j}\in \Rbb^{m\times m}$ and $p_{j,k}\in\Rbb$, with $k=0\ldots n$, be the coefficients of the trigonometric polynomials $\Qb_j$ and $\pb_j$, respectively. Notice that
 \al{\tilde J_\lambda(\pb_j,\Qb_j)=\int_\Tbb& \pb_j D_{IS}(\pb_j^{-1}\Qb_j,I)-\tr(\Qb_j )\nn\\
 &+\sum_{k=0}^n\tr( Q_{j,k}^T R_k)-\sum_{\substack{k=0}}^n p_{j,k} c_k+g_\lambda(\pb_j)\nn} 
  where \al{\label{ItaSai}D_{IS}(X,Y)=\log\det( X^{-1}Y)+\tr(XY^{-1})-m} denotes the Itakura-Saito distance between $X,Y\in\Sc^m_+$. Since $\pb_j(e^{i\theta})> 0$ and $D_{IS}(\pb_j^{-1}\Qb_j)\geq 0$ on $\Tbb$, it follows that 
 \al{\tilde J_\lambda & (\pb_j,\Qb_j) \geq\nn\\ &-\int_\Tbb\tr(\Qb_j )+\sum_{k=0}^n\tr( Q_{j,k}^T R_k)-\sum_{\substack{k=0}}^n p_{j,k} c_k+g_\lambda(\pb_j)\nn\\ & \rightarrow \infty\nn} 
 where we exploited the fact that $g_\lambda(\pb_j)\rightarrow \infty $ as $j\rightarrow\infty$, while the other terms converge to a finite value. \\
 {\em  Second case: $\pb\in \Bk^o_+$ and $\Qb\in \partial \Bk^m_+(E)$.} In this case the limit of $\tilde J_\lambda(\pb_j,\Qb_j)$ takes a finite value, however we will show that $(\pb,\Qb)$ cannot be a point of minimum for $\tilde J_\lambda$. By continuity we have 
 \al{\pb(e^{i\theta})\geq \mu:= \min_{\theta\in\Tbb} \pb(e^{i\theta})> 0.\nn}
 Now, consider the first variation of $\tilde J_\lambda(\pb,\cdot )$ at $\Qb+tI\in\Bk_+^m(E)$, with $t>0$, along the direction $\delta \Qb=I$: 
 \al{\delta \tilde J_\lambda(\pb,\Qb+tI; I)&=-\int_\Tbb \pb \tr(\Qb+tI)^{-1}+\tr(R_0)\nn\\
 &\leq -\mu \int_\Tbb \tr(\Qb+tI)^{-1}+\tr(R_0).\nn }
 By the Lebesgue's monotone convergence theorem, we have 
 \al{\underset{ t\rightarrow 0^+}{\lim}\delta \tilde J_\lambda(\pb,\Qb+tI; I) \leq  -\mu \int_\Tbb \tr (\Qb^{-1})+\tr(R_0)=-\infty\nn}
 because $\mu>0$ and $\tr(\Qb^{-1})$ is a positive rational function having some poles on $\Tbb$. We conclude that $(\pb,\Qb)$ cannot be a point of minimum for $\tilde J_\lambda$. 
   \qed\\
 
 We are now ready to prove the existence result.  \smallskip\\
{\em Proof of Theorem \ref{teo_Demp_gen}:} Consider the nonempty sublevel set 
{\small \al{\tilde J^{-1}_\lambda (-\infty,r] :=\{ (\pb,\Qb)\in\overline{\Bk^o_+}\times \overline{\Bk^m_+}(E) \hbox{ s.t. } \tilde J_\lambda(\pb,\Qb)\leq r\}\nn} }where the real number $r$ is taken sufficiently large. Such sublevel set is closed and bounded by Proposition \ref{prop4} and Proposition \ref{prop3}, respectively, and thus compact because $(\pb,\Qb)$ belongs to a finite dimensional vector space. Since $\tilde J_\lambda$ is lower semicontinuous, see Proposition \ref{prop2}, it follows that the regularized dual problem (\ref{dual_pb_r}) admits solution by the Weierstrass' theorem. Moreover, the minimum belongs to the interior of the feasible set  $\tilde J^{-1}_\lambda (-\infty,r]\subset \Bk^o_+\times \Bk_+^m(E)$ by Proposition \ref{prop4} and \ref{prop3}. Finally, this point of minimum is unique because $\tilde J_\lambda$ is strictly convex over $\Bk^o_+\times \Bk^m_+(E)$ by Proposition \ref{pro1}. Moreover, this point of minimum satisfies the stationarity conditions (i.e the first  variation of $\tilde J_\lambda$ must vanish along any direction) which implies the moment conditions in (\ref{optcond}). \qed\\

\section{Connection with the ML estimator}\label{sec:ML}
Consider a graphical model $\Gc(V,E)$ associated to the  ARMA process $y$
whose spectral density $\Phi\in\Sbb_+^m$ is of the form in (\ref{opt_form}). 
 Assume that  the data $y^N:=\{y(1)\ldots y(N)\}$ generated by this model is available. Then, an estimate of this graphical model can be obtained through Problem (\ref{moment_pb}): it is sufficient to replace $R_k$s and $c_k$s with their sample estimators using $y^N$
 \al{\hat R_k&=\hat R_{-k}^T=\frac{1}{N}\sum_{t=1}^{N-k} y(t+k)y(t)^T,\nn\\ \hat c_k&=\hat c_{-k}=\int_\Ts e^{i\theta k}\log\det(\hat \Phi_P )\nn}
 where 
 \al{\hat\Phi_P (e^{i\theta})=\sum_{k=1-h(N)}^{h(N)-1} \hat R_k e^{-i\theta k}\nn}
 denotes the $f$-periodogram computed from the data $y^N$ with $h(N)$ such that $h(N)\rightarrow \infty$ and $h(N)^2/N\rightarrow 0$ as $N\rightarrow \infty$. It is worth noting that $\hat \Phi_P$ is a mean-square consistent estimator under mild assumptions on $\Phi$,  \cite{FALCONI2023110672}. The dual problem coincides with (\ref{dual_pb}) where $R_k$s and $c_k$s are replaced with $\hat R_k$s and $\hat c_k$s. Let 
\al{\hat\Phi_n (e^{i\theta})&:=\sum_{k=-n}^{n} \hat R_k e^{-i\theta k}\nn\\
\hat \psi(e^{i\theta})&:=\sum_{k=-h(N)+1}^{h(N)-1} \hat c_k e^{-i \theta k}\approx \log\det(\hat\Phi_P).\nn}
It is not difficult to see that the dual function in (\ref{dual_pb}) can be written as  
\al{J(\pb,\Qb)=\int_\Ts -\pb \log\det(\pb^{-1}\Qb)+\tr(\Qb\hat \Phi_P)-\pb \hat \psi-\pb m\nn}
 where we exploited the identities 
 \al{&\tr\sum_{k=0}^n Q_k^T\hat R_k=\tr \int_\Ts \Qb \hat \Phi_n=\tr\int_\Ts \Qb \hat 	\Phi_P.\nn\\
&\tr\sum_{k=0}^n p_k \hat c_k=\int_\Ts \pb \sum_{k=-n}^n  \hat c_k e^{-i \theta k}=\int_\Ts \pb \hat \psi\nn } 
because  $\pb\in\overline{\Bk_+^o}$ and $\Qb\in\overline{\Bk_+^m}(E)$. Thus, 
\al{J&(\pb,\Qb)=\int_\Ts \pb \log\det(\pb\Qb^{-1})+\tr(\Qb\hat \Phi_P)-\pb \hat \psi-\pb m\nn\\
&\approx \int_\Ts \pb [\log\det(\pb\Qb^{-1})+\tr(\pb^{-1}\Qb\hat \Phi_P)-\log \det(\hat\Phi_P)-m] \nn\\
&=\int_\Ts \pb D_{IS}(\hat\Phi_P,\pb\Qb^{-1})\nn}
where we recall that $D_{IS}$ is the Itakura-Saito distance  defined in (\ref{ItaSai}). Let $(\pb,\Qb)\in \Bk^o_+ \times \Bk^m_+(E)$. 
Then, $\pb(e^{i\theta})>0$ and $D_{IS}(\hat\Phi_P,\pb\Qb^{-1})\geq 0$ on $\Tbb$. By the weighted mean value theorem it follows that there exists $\bar \theta$ such that 
\al{\label{ineq_J}J&(\pb,\Qb)\approx \int_\Ts \pb D_{IS}(\hat\Phi_P,\pb\Qb^{-1})\nn\\ &=\pb(e^{i\bar \theta})\int_\Ts   D_{IS}(\hat\Phi_P,\pb\Qb^{-1})\\&=\pb(e^{i\bar \theta})[\int_\Ts \log\det(\pb\Qb^{-1})+\tr(\pb^{-1}\Qb\hat \Phi_P)\nn\\
& \hspace{3.8cm}-\log \det(\hat\Phi_P)-m].\nn}

\prop Consider an ARMA process with spectral density $\Phi=\pb\Qb^{-1}$ with $\pb\in\Bk_+^o$ and $\Qb\in\Bk_+^m(E)$.
For $N$ sufficiently large, the negative log-likelihood of the data $y^N$ under this ARMA model  is (up to constant terms not depending on $\pb$ and $\Qb$)
\al{\label{lik_model}\ell(y^N;\pb,\Qb)\approx \frac{N}{2}\int_\Ts \log\det(\pb\Qb^{-1})+\tr(\pb^{-1}\Qb\hat \Phi_P). } 
\eprop

\IEEEproof  Let $\Wb$ be a causal rational right spectral factor of $\pb^{-1}\Qb$ which is analytic on $\Tbb$. Let $s$ be the practical length of $\Wb$, that is a natural number for which 
\al{\label{high_ord}\Wb(z)\approx\sum_{k=0}^s W_kz^{-k}.}
Since $\pb(z) \Qb(z)^{-1}=\Wb^{-1}(z)\Wb(z^{-1})^{-T}$, then a model for the ARMA process is
\al{\label{mod1}\Wb(z)y(t)=e(t)}
where $e$ is normalized white Gaussian noise. By (\ref{high_ord}), we can approximate (\ref{mod1}) with the high order AR process
\al{\sum_{k=0}^sW_k y(t-k)=e(t).\nn} 
The negative log-likelihood of this model is, \cite{whittle1953estimation},
\al{\ell(&y^N;\pb,\Qb)\nn\\
&= \frac{N-s}{2}\int_\Ts -\log \det(\Gb^*\Gb)+\frac{N}{N-s}\tr(\Gb^*\Gb\hat \Phi_P)\nn\\
&\approx \frac{N}{2}\int_\Ts -\log \det(\Gb^*\Gb)+\tr(\Gb^*\Gb\hat \Phi_P)\nn}
where $\Gb(z)=\sum_{k=0}^sW_kz^{-k}$ and we exploited the fact that $N\gg s$. Finally, since $\Gb\approx \Wb$ and thus $\Gb^*\Gb\approx \Wb^*\Wb=\pb^{-1}\Qb$ we obtain (\ref{lik_model}).
\qed\\

Taking into account (\ref{ineq_J}), we conclude that for $N$ sufficiently large
\al{\label{ineq_lik_J}\ell (y^N;\pb,\Qb)\approx 	 \frac{N}{2} \left[\frac{J(\pb,\Qb)}{\pb(e^{i\bar\theta})}+\int\log\det(\hat \Phi_P)+m\right].} Although Problem (\ref{moment_pb}) provides a tool to estimate an ARMA graphical model from the data, it is required to know the topology, i.e. the set $E$. In the following we will show how to address the problem in the case the topology is not known. 

 \section{A Bayesian viewpoint}\label{sec:Bayes}
 The aim of this section is to show how to regularize the dual problem (\ref{dual_pb})
 in order to estimate an ARMA graphical model in the case $E$ is not known. More precisely we take the Bayesian perspective proposed in \cite{ZORZI2019108516} and thus we model $\pb$ and $\Qb$ as a stochastic trigonometric polynomials (i.e. $p_k$s and the entries of $Q_k$s
 are random variables) taking values in $\overline{\Bk_+^o}$ and  $\overline{\Bk_+^m}:=\overline{\Bk_+^m}(V\times V)$, respectively. We also assume they are independent. Let $f(\Qb|\bm \gamma)$ be the prior of $\Qb$ where $\bm \gamma$ denotes the vector containing $\gamma_{jh}>0$, $j\geq h$, which are referred to as hyperparameters.
 Recall that $[\Qb]_{jh}$ is the entry of $\Qb$ in position $(j,h)$, that is 
 $[\Qb(e^{i\theta})]_{jh}=[Q_0]_{jh}+0.5\sum_{k=1}^n [Q_k]_{jh}e^{-i\theta k}+[Q_k]_{hj}e^{i\theta k}$. We assume that $[\Qb]_{jh}$s are independent each other. Thus, $f(\Qb|{\bm \gamma})=\prod_{j\geq h} f([\Qb]_{jh}|{\bm \gamma})$ and $f([\Qb]_{jh}|{\bm \gamma})=e^{-\gamma_{jh}q_{jh}(\Qb)} /c_{jh}$ where 
  \al{q_{jh}(\Qb):=\max \{|[Q_0]_{jh}|,\max_{k=1\ldots n} |[Q_k]_{jh}|,\max_{k=1\ldots n} |[Q_k]_{hj}|\}\nn}  and
$c_{jh}$ is the normalizing constant. Therefore, we obtain 
\al{f(\Qb|{\bm \gamma})= e^{-h_W(\Qb)}/\prod_{j\geq h} c_{jh}\nn}
where \al{h_W(\Qb)=\sum_{j\geq h} \gamma_{jh}q_{jh}(\Qb).\nn}
Moreover, we model $\gamma_{jh}$s as independent random variables with exponential distribution:
\al{f(\bm \gamma)=\prod_{j\geq h} \varepsilon e^{-\varepsilon \gamma_{jh}}\nn}
where $\varepsilon>0$ is a fixed small constant.
Let $f(\pb)$ be the prior of $\pb$. We assume that 
\al{\label{priorP}f(\pb)= c^{-1} e^{-g_\lambda(\pb) }}
where $c$ is the normalizing constant and $\lambda>0$ is the corresponding hyperparameter. It is worth noting that $\lambda$ plays the same role of the regularization parameter in (\ref{dual_pb_r}) and thus we assume it is fixed. 

\begin{example}
Consider the trigonometric polynomial $\pb(e^{i\theta})=1+p_1cos(\theta)+p_2\cos(2\theta)$ with $n=2$. Figure \ref{fig:priorP} shows the prior (\ref{priorP}) for $\lambda=0.2$ and $\lambda=1$. We notice that the larger $\lambda$ is, the more the neighborhood of high-probability around $\pb(e^{i\theta})=1$ is shrunk. \begin{figure}
\centering
\includegraphics[width=0.5\textwidth]{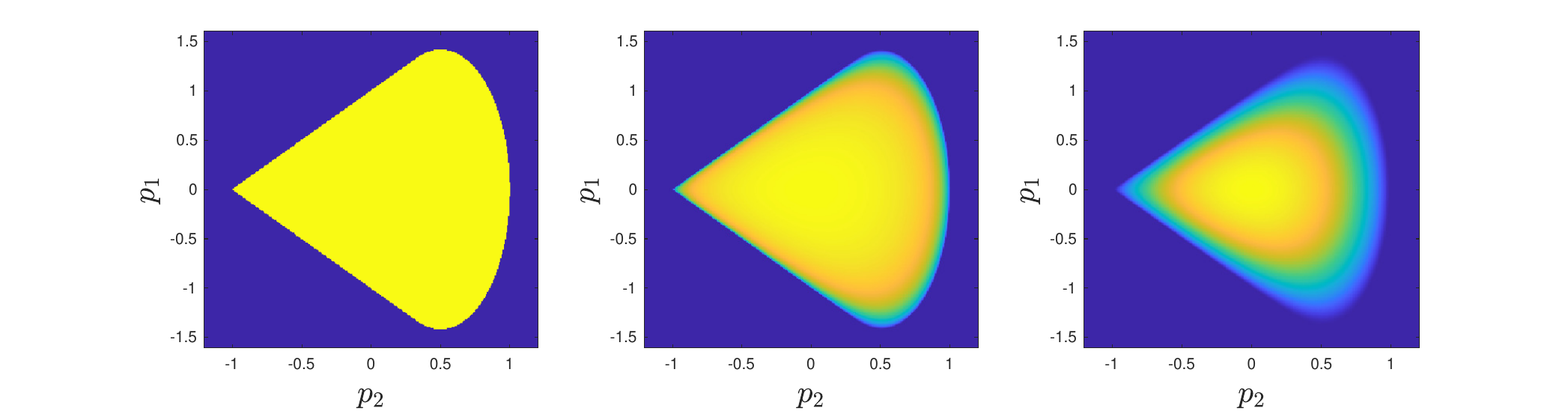}
\caption{{\em Left.} The yellow region denotes the values of $p_1$ and $p_2$ such that the corresponding $p\in \Bk_+^o$;   {\em Center.} Prior  with $\lambda=0.2$,  hot colors indicate high probability while cold colors low probability; {\em Right.} Prior  with $\lambda=1$.}\label{fig:priorP}
\end{figure}

\end{example}
The negative log-likelihood of $y^N$, $\pb$, $\Qb$ and $\bm \gamma$ takes the form
\al{\ell&(y^N, \pb,\Qb, \bm\gamma)=-\log f(y^N,\pb,\Qb,{\bm \gamma})\nn\\
&=-\log f(y^N|\pb,\Qb)-\log f(\pb)-\log f(\Qb|{\bm \gamma})-\log f(\bm \gamma).\nn} 
Since, $-\log f(y^N|\pb,\Qb)$ coincides with $\ell(y^N; \pb,\Qb)$ (i.e. the likelihood function considered in Section \ref{sec:ML}), we have 
\al{\label{log_lik_Bayes}\ell(y^N, \pb,\Qb; \bm \gamma)=\ell & (y^N;\pb,\Qb)+g_\lambda(\pb)+h_W(\Qb)\nn\\ &+\sum_{j\geq h}\log c_{jh}+\varepsilon \gamma_{jh}}
where we have neglected the constant terms not depending on $\pb$, $\Qb$ and $\bm \gamma$. The maximum a posteriori  (MAP) estimator of $(\pb,\Qb)$ is given by minimizing $\ell(y^N, \pb,\Qb, \gamma)$ with respect to $(\pb,\Qb)\in\overline{\Bk_+^o}\times \overline{\Bk_+^m}$. The hyperparameter vector $\bm \gamma$, however, is not known.  
According to the empirical Bayes method  \cite{RASMUSSEN_WILLIAMNS_2006}, it is possible to find an estimate  $\hat{\bm \gamma}$ by minimizing the so called negative log-marginal likelihood with respect to $\bm \gamma$:
\al{\label{marginal_lik}-\log \int_{\overline{\Bk_+^o}\times \overline{\Bk_+^m}} f(y^N,\pb,\Qb,{\bm \gamma})\mathrm d \pb \mathrm d \Qb;} then, the MAP estimator of $(\pb,\Qb)$ is given by minimizing $\ell(y^N,\pb,\Qb,\hat{\bm \gamma})$. However, it is difficult to find an analytic expression for (\ref{marginal_lik}) making challenging the optimization of $\bm \gamma$.  One way to overcome this issue is to consider an upper bound for (\ref{log_lik_Bayes}). By (\ref{ineq_lik_J}) and \cite[Proposition 3]{ZORZI2019108516}, the following inequality approximately holds for $N$ sufficiently large
\al{\label{ineq_l_lt}\ell(y^N, \pb,\Qb, \bm\gamma)\leq \tilde \ell(y^N, \pb,\Qb,\alpha, \bm\gamma)}
where 
\al{\label{upper_bound}\tilde \ell( y^N, \pb,\Qb,\alpha, &\bm\gamma):=\frac{N}{2 }[\alpha^{-1}J(\pb,\Qb)+\int\log\det(\hat \Phi_P)+m ] \nn\\ 
&+g_\lambda(\pb)+h_W(\Qb) -\sum_{j>h}(2n+1)\log\gamma_{jh}\nn\\ &-\sum_{j=1}^n (n+1)\log \gamma_{jj}+\sum_{j\geq h}\varepsilon \gamma_{jh},}
$\alpha= \pb(e^{i\bar \theta})$ and we neglected the constant terms in the inequality (\ref{ineq_l_lt}). Then, a simplified approach to estimate $\bm \gamma$ is the generalized maximum likelihood (GML) method \cite{zhou1997approximate}: this hyperparameter is jointly optimized with $\pb$ and  $\Qb$ minimizing (\ref{upper_bound}). However, the optimization with respect to $\pb$ is difficult because of the presence of $\alpha$. Accordingly, we can find a preliminary estimate $(\pb^{(0)},\Qb^{(0)})$ and then, in view of (\ref{ineq_J}), an estimate of $\alpha$ is given by

{\small \al{\label{alpha1}\hat\alpha&=\frac{J(\pb^{(0)},\Qb^{(0)})}{\int_\Ts \log\det(\pb^{(0)}(\Qb^{(0)}\hat\Phi_P)^{-1})+\tr((\pb^{(0)})^{-1}\Qb^{(0)}\hat \Phi_P)-m }.}} 

\teo \label{teo_nuovo2} Let $\alpha>0$ and assume that there exists a power spectral density $\bar\Phi$ which is positive definite on some open interval in $\Tbb$ and  such that
 \al{\label{esist2}\int_\Ts \bar\Phi e^{i\theta k}=R_k, \quad k=0\ldots n. } Then, the  optimization problem
\al{\underset{(\pb,\Qb,\bm \gamma)\in\mathcal C }{\min}\,\tilde \ell( y^N, \pb,\Qb,\alpha, &\bm\gamma),\nn}
with $\mathcal C:=\{(\pb,\Qb,\bm \gamma) \hbox{ s.t. } 	\pb \in\overline{\Bk_+^o}, \; \Qb\in\overline{\Bk_+^m},\;  \gamma_{jh}> 0\}$, admits solution. Moreover, the latter is such that $(\pb,\Qb)\in \Bk_+^o\times \Bk_+^m$. \eteo

The proof is deferred after Theorem \ref{teo_nuovo}. Since the joint optimization of $\pb,\Qb,\bm\gamma$ is still difficult, we perform the optimization of $(\pb,\Qb)$ and $\bm\gamma$ iteratively through a two-step algorithm  where at the $l$-iteration the variables are estimated as follows: 
\al{\label{P1}(\hat \pb^{(l)},\hat\Qb^{(l)})&=\underset{{\small (\pb,\Qb)\in \overline{\Bk_+^o}\times \overline{\Bk_+^m}} }{\mathrm{argmin}}\,\tilde \ell(y^N, \pb,\Qb,\hat\alpha, \hat{\bm\gamma}^{(l)})\\
\label{P2}\hat{\bm\gamma}^{(l+1)}&=\underset{\gamma_{jh>0}}{\mathrm{argmin}}\,\tilde \ell(y^N,\hat \pb^{(l)},\hat \Qb^{(l)},\hat\alpha,  \bm\gamma ).}
It is not difficult to see that Problem (\ref{P1}) is equivalent to the following problem 
\al{\label{dual_pb_new}& \underset{\pb,\Qb}{ \min\;} J(\pb,\Qb)+\frac{2\alpha}{N}[g_\lambda(\pb)+h_W(\Qb) ]\\ & \hbox{ s.t. }  \pb\in \overline{\Bk_+^o},\; \Qb\in\overline{\Bk_+^m} \nn 
  }
with $\alpha=\hat \alpha$ and $\bm\gamma=\hat{\bm\gamma}^{(l)}$ which is a regularized version of the dual problem in (\ref{dual_pb}). In particular: the regularization term $g_\lambda(\pb)$ avoids that the optimal solution $\pb$ is on the boundary, while $h_W(\Qb)$ induces sparsity in $\Qb$ in a soft way; finally, the factor $2\alpha/N$ tunes the best trade-off between the ``fit function'' and the regularization terms. 

\teo\label{teo_nuovo} Let $\alpha>0$ and $\bm \gamma$ be such that $\gamma_{jh}\geq 0$. Moreover, assume that there exists a power spectral density $\bar\Phi$ which is positive definite on some open interval in $\Tbb$ and  such that (\ref{esist2}) holds.
 Then, Problem~(\ref{dual_pb_new}) admits a unique solution $(\hat \pb,\hat \Qb)\in \Bk_+^o\times \Bk_+^m$.\eteo
\IEEEproof First, notice that the results in Section \ref{sec_dualpb} also holds in the particular case $E=V\times V$. Moreover, the objective function in Problem (\ref{dual_pb_new}) can be written as 
\al{\check J_{\lambda,\alpha,\bm \gamma}(\pb,\Qb)=\tilde J_\frac{2\alpha\lambda}{N}(\pb,\Qb)+\frac{2\alpha}{N}h_W(\Qb)\nn}
where $\tilde J$ is the regularized dual function defined in (\ref{dualfreg}) and recall that $\lambda>0$.
In plain words, the difference between $\tilde J$ and $\check J $ is that in the latter appears the term $s(\Qb):=2\alpha N^{-1}h_W(\Qb)$. Then, it follows that $\check J$ is strictly convex in  $\Bk_+^o\times \Bk_+^m$ because $\tilde J$ is strictly convex, by Proposition \ref{pro1}, and $s(\Qb)$ is convex. Moreover, $\check J$ is lower semicontinuous on $\overline{\Bk_+^o}\times \overline{\Bk_+^m}$ because $\tilde J$ is lower semicontinuous, by Proposition \ref{prop2}, and $s(\Qb)$ is continuous. 

Next, take a sequence  $\{(\pb_j,\Qb_j)\in\overline{\Bk^o_+}\times \overline{\Bk^m_+}, \, j\in\Nbb\}$ such that $\|(\pb_j,\Qb_j)\|\rightarrow \infty$ as $j\rightarrow\infty$. Then, in view of Proposition \ref{prop3}, we have 
\al{\underset{j\rightarrow \infty }{\mathrm{lim\,inf}} \check J_{\lambda,\alpha,\bm \gamma}(\pb_j,\Qb_j)\geq \underset{j\rightarrow \infty }{\mathrm{lim\,inf}} \tilde J_\frac{2\alpha\lambda}{N}(\pb_j,\Qb_j)\rightarrow \infty.\nn} Thus, this sequence cannot be an infimizing one.

Finally, consider a sequence $\{(\pb_j,\Qb_j)\in {\Bk^o_+}\times {\Bk^m_+}, \, j\in\Nbb\}$ converging to $(\pb,\Qb)$ which belongs to the boundary of $\overline{\Bk^o_+}\times \overline{\Bk^m_+}$. The latter cannot be an infimizing sequence. This statement can be proved using the same arguments exploited in the proof of Proposition \ref{prop4}: the unique substantial difference regards the case $\pb\in \Bk_+^o$ and $\Qb\in\partial \Bk_+^m$. The function $\check J_{\lambda,\alpha,\bm \gamma}(\pb,\cdot)$ is differentiable at $\Qb+tI$, with $t>0$, along the direction $\delta \Qb=I$. More precisely, the first variation is 
\al{\delta \check J_{\lambda,\alpha,\bm \gamma}(&\pb,\Qb+tI; I)\nn\\ &=-\int_\Tbb \pb \tr(\Qb+tI)^{-1}+\tr(R_0)+\sum_{j} \frac{2\alpha\gamma_{jj}}{N}\nn\\
 &\leq -\mu \int_\Tbb \tr(\Qb+tI)^{-1}+\tr(R_0)+\sum_{j} \frac{2\alpha\gamma_{jj}}{N}\nn}
and by the Lebesgue's monotone convergence theorem, we have 
 \al{\underset{ t\rightarrow 0^+}{\lim}&\delta \check J_{\lambda,\alpha,\bm \gamma}(\pb,\Qb+tI; I)\nn\\ &\leq  -\mu \int_\Tbb \tr (\Qb^{-1})+\tr(R_0)+\sum_{j} \frac{2\alpha\gamma_{jj}}{N} =-\infty.\nn}

In the light of the above results, taken $r\in \Rbb$ sufficiently large, we can restrict the search of the minimum over the sublevel set  $\{\, (\pb,\Qb)\in  \Bk^o_+ \times \Bk^m_+  \hbox{ s.t. } \check J_{\lambda,\alpha,\bm \gamma}(\pb,\Qb)\leq r\,\}$ which is compact. Since $\check J$ is strictly convex over this sublevel set and by the Weierstrass' theorem, we conclude that Problem (\ref{dual_pb_new}) admits solution which is also unique.
\qed\\

{\em Proof of Theorem \ref{teo_nuovo2}:}  Since $\mathcal C$ is open and unbounded, we need to show that we can restrict the search of the minimum of $\tilde \ell$ over a compact sect $\mathcal C^\star\subset \mathcal C$ for which $\tilde \ell$ is continuous, then the claim follows by the Weierstrass' theorem. Moreover, we must show that $(\pb,\Qb)\in \Bk_+^o\times \Bk_+^m$ for any $(\pb,\Qb,\bm \gamma)\in{\cal C}^\star$. Notice that 
\al{\label{ineqell}\tilde \ell(& y^N, \pb,\Qb,\alpha, \bm\gamma)=\frac{N}{2\alpha } \check J_{\lambda,\alpha,\bm \gamma}  (\pb,\Qb)+\kappa   \nn\\ 
& -\sum_{j>h}(2n+1)\log\gamma_{jh}-\sum_{j=1}^n (n+1)\log \gamma_{jj}+\sum_{j\geq h}\varepsilon \gamma_{jh}\nn\\
&\hspace{-0,5cm}\geq \frac{N}{2\alpha } \check J_{\lambda,\alpha,\bm 0}  (\pb,\Qb)+\kappa -\sum_{j>h}(2n+1)\log\gamma_{jh}  \nn\\ &-\sum_{j=1}^n (n+1)\log \gamma_{jj}+\sum_{j\geq h}\varepsilon \gamma_{jh}} 
where $\kappa$ contains the terms not depending on $\pb$, $\Qb$ and $\bm \gamma$. Now, consider a sequence  $\{(\pb_j,\Qb_j,\bm\gamma_j)\in\mathcal C, \, j\in\Nbb\}$ such that there exists at least one $\gamma_{jh}\rightarrow \infty$. Then, in view of (\ref{ineqell}),
\al{\lim_{j\rightarrow \infty}\tilde \ell(& y^N, \pb_j,\Qb_j,\alpha, \bm\gamma_j)=\infty\nn }
because, by Theorem \ref{teo_nuovo}, $\check J_{\lambda,\alpha,\bm 0}$ is bounded from below on $\overline{\Bk^o_+}\times \overline{\Bk^m_+}$.
Next, take a sequence  $\{(\pb_j,\Qb_j,\bm\gamma_j)\in {\mathcal C}, \, j\in\Nbb\}$ such that there exists at least one $\gamma_{jh}\rightarrow 0$. 
Then, by (\ref{ineqell}), we have 
that $\tilde \ell( y^N, \pb_j,\Qb_j,\alpha, \bm\gamma_j)\rightarrow\infty$ as $j\rightarrow \infty$  because $\check J_{\lambda,\alpha,\bm 0}$ is bounded from below. Therefore, we can restrict the search of the minimum over the set $\tilde{ \cal C}:=\{(\pb,\Qb,\bm \gamma) \in\mathcal C \hbox{ s.t. } \delta \leq\gamma_{jh}\leq \mu\}$ where $\mu>\delta>0$. Finally, by Theorem \ref{teo_nuovo}, we can restrict the search of the minimum over 
\al{\mathcal C^\star:=\{(\pb,\Qb,\bm \gamma) \in\mathcal C \hbox{ s.t. } (\pb,\Qb)\in{\cal A},\;\delta\leq \gamma_{jh}\leq \mu\}\nn}
and ${\cal A}\subset {\Bk^o_+}\times  {\Bk^m_+}$ is a compact set. Accordingly, $\mathcal C^\star$ is a compact set. \qed
\smallskip\\
It is not difficult to see that (\ref{P2}) is the same optimization problem in \cite[Section 4]{ZORZI2019108516} and it admits the following closed form solution:
\al{\label{formula_gamma}\hat \gamma_{jh}^{(l+1)}=\left\{\begin{array}{cc}\frac{n+1}{q_{jh}(\hat \Qb^{(l)})+\varepsilon }, & \hbox{ if $j=h$} \\ \frac{2n+1}{q_{jh}(\hat \Qb^{(l)})+\varepsilon }, & \hbox{ if $j>h$.}  \end{array}\right.
} The resulting estimator, hereafter called Generalized Maximum Likelihood (GML) estimator, is described in Algorithm~\ref{Algo1} where $\epsilon>0$ is a fixed threshold.  

\begin{algorithm}
    \caption{\textbf{GML} estimator}\label{Algo1}
     Compute  $(\pb^{(0)},\Qb^{(0)})$ as solution of (\ref{dual_pb_r}) with $E=V\times V$ 
    \newline
    Compute $\hat\alpha$ as in (\ref{alpha1})
    \newline
    Set $\hat{\bm \gamma}^{(1)}=0$ and $l=0$ \newline    
    \textbf{Repeat} 
    \newline
               \hspace*{0.4cm} $l=l+1$
            \newline
   \hspace*{0.4cm} Compute  $(\hat \pb^{(l)},\hat \Qb^{(l)})$ as  solution of (\ref{dual_pb_new})
        \newline
        \hspace*{0.4cm} Compute  $\hat{\bm\gamma}^{(l+1)}$ through  (\ref{formula_gamma})
            \newline
   \textbf{until} $\| \hat \Qb^{(l)}-\hat \Qb^{(l-1)}\|\leq \epsilon$ 

\end{algorithm}

\begin{corollary} \label{cor_opt}Let $(\hat \pb^{(l)},\hat \Qb^{(l)},\hat{\bm\gamma}^{(l+1)})$, with $l\in\mathbb N$, be the sequence generated by Algorithm \ref{Algo1} where $\hat \alpha$ is fixed.  If there is a limit point $(\hat \pb,\hat \Qb,\hat{\bm\gamma})$ of such a sequence, then it is a coordinatewise minimum point of $\tilde \ell$, i.e.
\al{\tilde \ell(y^N,\hat  \pb,\hat \Qb,\hat \alpha, \hat{\bm\gamma}) \leq \tilde \ell(y^N, \pb, \Qb,\hat \alpha, \hat{\bm\gamma})\nn\\
\tilde \ell(y^N,\hat  \pb,\hat \Qb,\hat \alpha, \hat{\bm\gamma}) \leq \tilde \ell(y^N, \hat \pb,\hat \Qb,\hat \alpha, \bm\gamma)\nn}
for any $\pb\in\overline{\Bk_+^o}$, $\Qb\in\overline{\Bk_+^m}$ and $\gamma_{jh}>0$.
\end{corollary}
\IEEEproof The statement is supported by the fact that at every stage of the sequential procedure, we are able to identify a specific variable for which the function $\tilde \ell$ attains its minimum value.
\qed

It is worth noting one could consider the following upper bound for $\ell(y^N, \pb,\Qb, \bm \gamma)$ which considers the likelihood of the data:
\al{\label{def_l_check} \check \ell  &(y^N, \pb,\Qb, \bm\gamma):=\ell(y^N;\pb,\Qb)+g_\lambda(\pb)+h_W(\Qb)\nn\\ & -\sum_{j>h} (2n+1)\log\gamma_{jh}-\sum_{j=1}^n (n+1)\log \gamma_{jj}+\sum_{j\geq h} \varepsilon\gamma_{jh}} 
with $\ell  (y^N; \pb,\Qb)$ given by the approximation in (\ref{lik_model}). This leads to a two-step algorithm: in the first step $\check \ell$ is optimized with respect to $(\pb,\Qb)$, while in the second step with respect to $\bm \gamma$. However, $\check \ell$ is not convex in $(\pb,\Qb)$ and thus it is plausible to obtain a local minimum which is not a global one in the first step. As a consequence, the sequence generated by the corresponding sequential algorithm does not enjoy the property in Corollary \ref{cor_opt}.
Accordingly, the upper bound $\tilde \ell$ aims to ``convexify'' the problem in $(\pb,\Qb)$ with the introduction of the additional parameter $\alpha$ and guarantees the sequence satisfies the property in Corollary \ref{cor_opt}.

\section{Numerical experiments}\label{sec:example}

\subsection{Local minima using $\check \ell$}
As we noticed in Section \ref{sec:Bayes}, an upper bound for the negative log-likelihood $\ell(y^N,  \pb,\Qb, \gamma)$ is given also by $\check \ell$ defined in (\ref{def_l_check}). However, the minimization of $\check \ell$ with respect to $(\pb,\Qb)$ could lead to a local minimum which is not a global one. Next, we provide a numerical example in which this situation happens.

\begin{figure}
\centering
\includegraphics[width=0.24\textwidth]{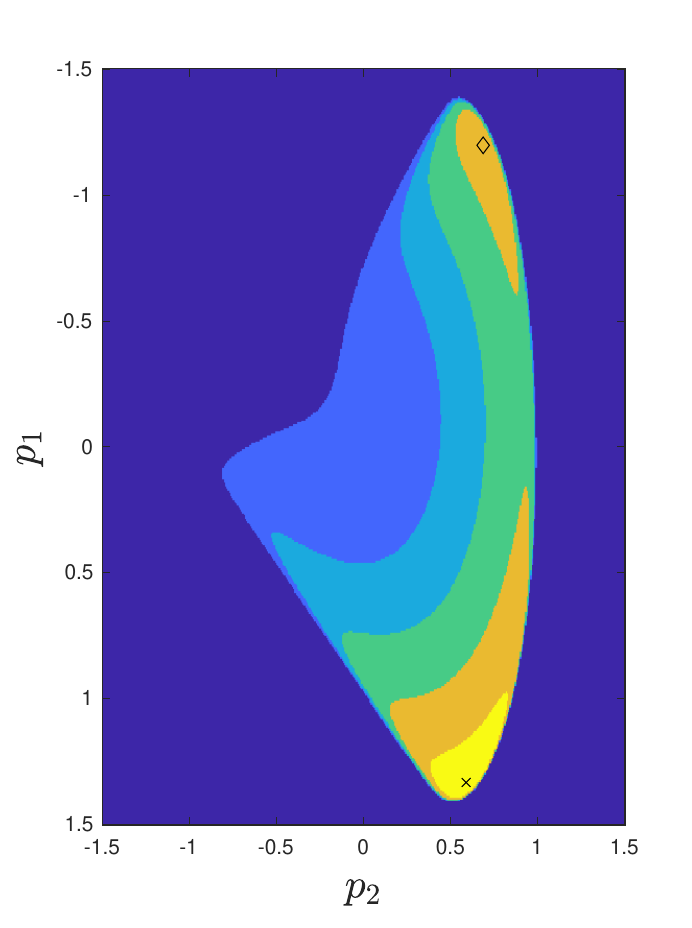}
\includegraphics[width=0.24\textwidth]{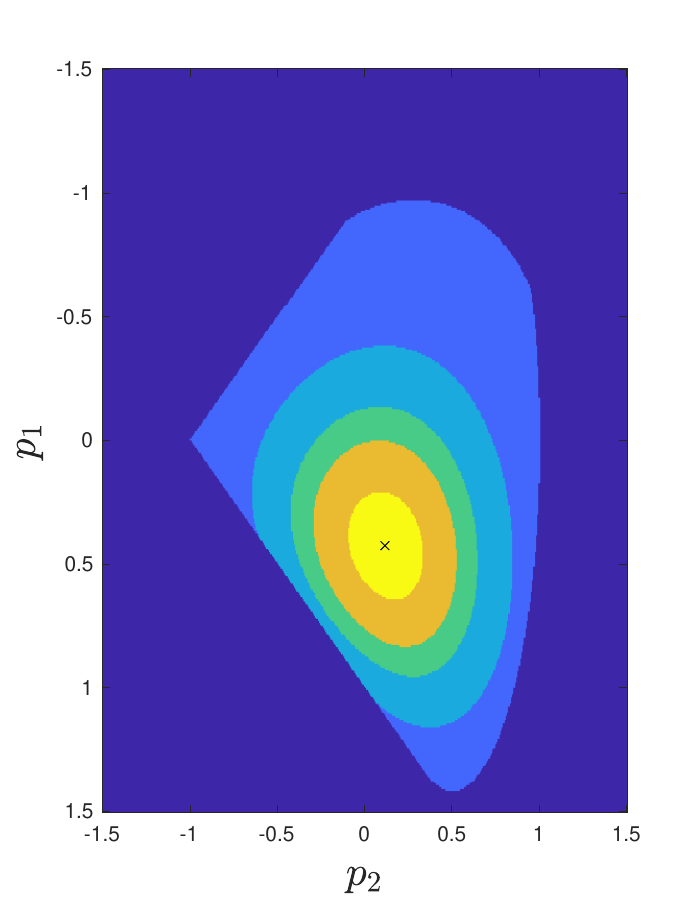}
\caption{Some sublevel sets of $(p_1,p_2) \mapsto \check\ell(y^N, \pb,\Qb^\star, \bm\gamma)$, left, and $(p_1,p_2)  \mapsto  \tilde \ell(y^N, \pb,\Qb^\star, \bm\gamma)$, right. 
The sublevel sets depicted with hot colors correspond to a smaller constant than the one of the sublevel sets  depicted with cold colors. The former are overlapped to the latter. The global minima are denoted by a cross, while the other local minima by a diamond.
}\label{fig:level}
\end{figure}

\begin{figure}
\centering
\includegraphics[width=0.5\textwidth]{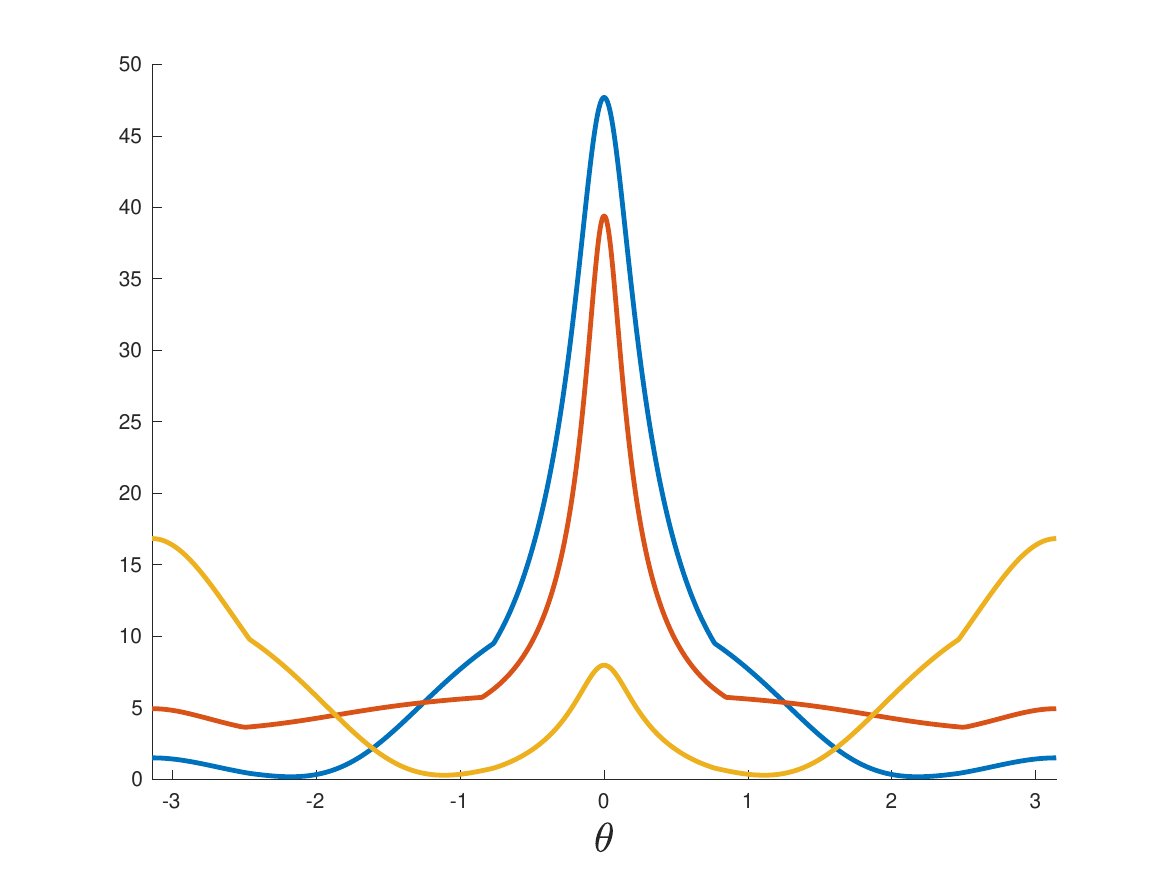} 
\caption{Pointwise Frobenius norm of $\Phi^\star=\pb(\Qb^\star)^{-1}$ when: $(\pb,\Qb^\star)$ is the global minimum of $\check \ell(y^N,\cdot,\cdot,\bm \gamma)$ (blue line); $(\pb,\Qb^\star)$ is  the other local minimum of $\check \ell(y^N,\cdot,\cdot,\bm \gamma)$ (yellow line); $(\pb,\Qb^\star)$ is the global minimum of $\tilde \ell(y^N,\cdot,\cdot,\hat\alpha,\bm \gamma)$ (red line). 
}\label{fig:level2}
\end{figure}

We consider an ARMA graphical model with $m=6$ nodes, oder $n=2$ and with $E:=\{(1,1),$ $(1,4),\,(1,6),\,(2,2),\,(3,3),\,(3,4),\,(3,6),\,(4,4),\, (4,6),\,(5,5),$ $(5,6),\,(6,6) \}$. From this model we generate a data sequence $y^N$ of length $N=500$. Then, we fix the hyperparameter vector $\bm\gamma$ as follows:
\al{\gamma_{jh}= \left\{\begin{array}{ll}10^3, & \hbox{ if }  (j,h)\in E  \\5\cdot 10^4, &     \hbox{ if }  (j,h)\notin E \end{array}\right.\nn}  
that is we consider a hyperparameter which induces the correct sparsity pattern on $\Qb$, but it induces  a little more bias than required in the other entries. Finally, we set  $\lambda=1$. Then, we solve numerically  the minimization problem 
\al{\label{numericalopt}(\hat \pb^\star,\hat\Qb^\star)&=\underset{(\pb,\Qb)\in\overline{\Bk_+^o}\times\overline{\Bk_+^m}}{\mathrm{argmin}}\,\check \ell(y^N, \pb,\Qb, \bm\gamma).} Figure \ref{fig:level}(left) shows some sublevel sets of the function $(p_1,p_2) \mapsto \check\ell(y^N, \pb,\Qb^\star, \bm\gamma)$, recall that $\pb$ is characterized by the variables $p_1,p_2$ see (\ref{defpQ}). We can see that there are two  local minima and one of them is not a global minimum. Now, we do the same thing with $\check \ell$. More precisely, we solve numerically  the minimization problem  (\ref{numericalopt}) where replace the objective function with $\tilde\ell(y^N, \pb,\Qb, \hat \alpha,\bm\gamma)$ and $\hat \alpha$ is computed as in the preliminary step of Algorithm \ref{Algo1}. Figure \ref{fig:level}(right) shows some sublevel sets of the function  $(p_1,p_2) \mapsto \tilde\ell(y^N, \pb,\Qb^\star, \hat \alpha,\bm\gamma)$. In this case we only have the  global minimum.   
Figure \ref{fig:level2} shows the pointwise Frobenius norm of the spectrum $\Phi^\star=\pb (\Qb^\star)^{-1}$ when: $(\pb,\Qb^\star)$ is the global minimum of $\check \ell$; $(\pb,\Qb^\star)$ is  the other local minimum of $\check \ell$; $(\pb,\Qb^\star)$ is the global minimum of $\tilde \ell$. While the two spectra corresponding to the global minima of $\check \ell$ and $\tilde \ell$ are similar, the spectrum corresponding to the local minimum of $\check \ell$ is very different. 

\subsection{Performance analysis}
\begin{figure}
\centering
\includegraphics[width=0.24\textwidth]{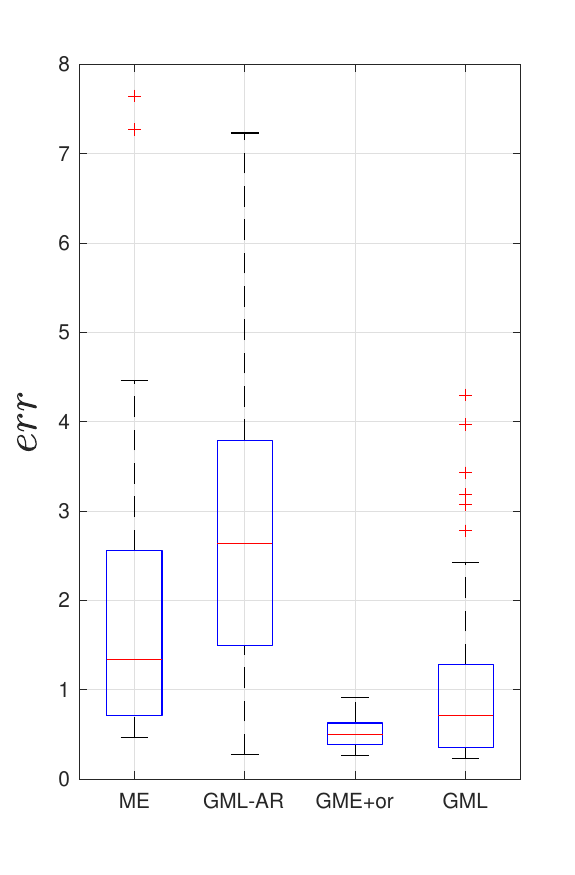}
\includegraphics[width=0.24\textwidth]{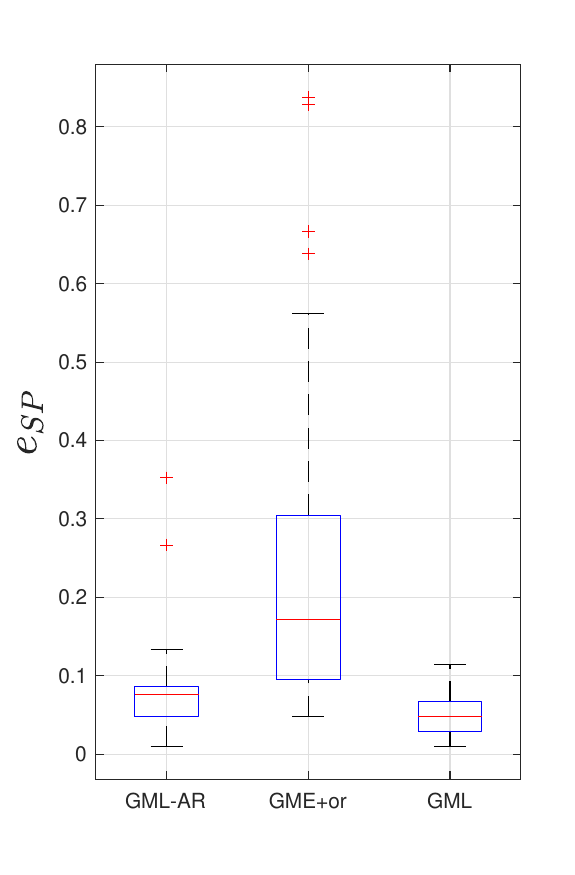}
\caption{First Monte Carlo experiment with $N=500$. Boxplots of the relative error (left) and the fraction of misspecified edges (right).}\label{figN500}
\end{figure}

\begin{figure}
\centering
\includegraphics[width=0.24\textwidth]{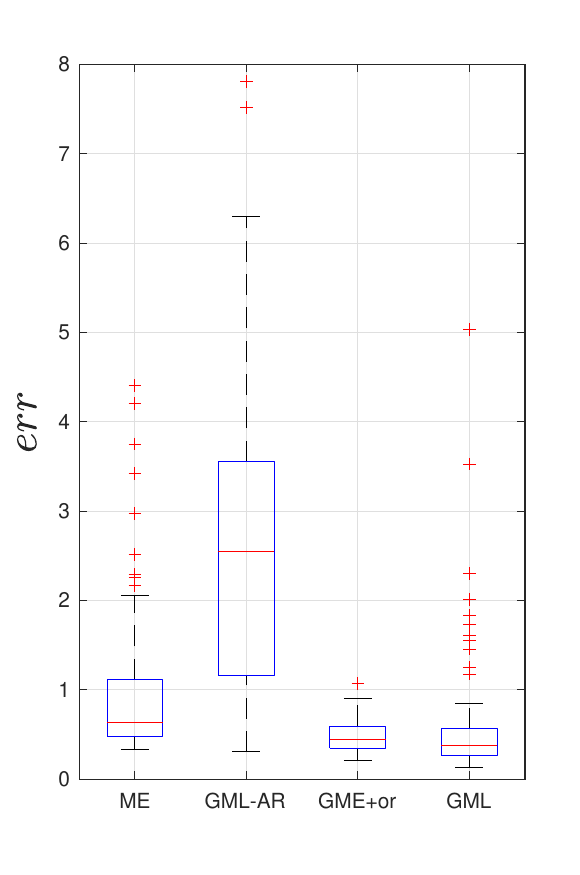}
\includegraphics[width=0.24\textwidth]{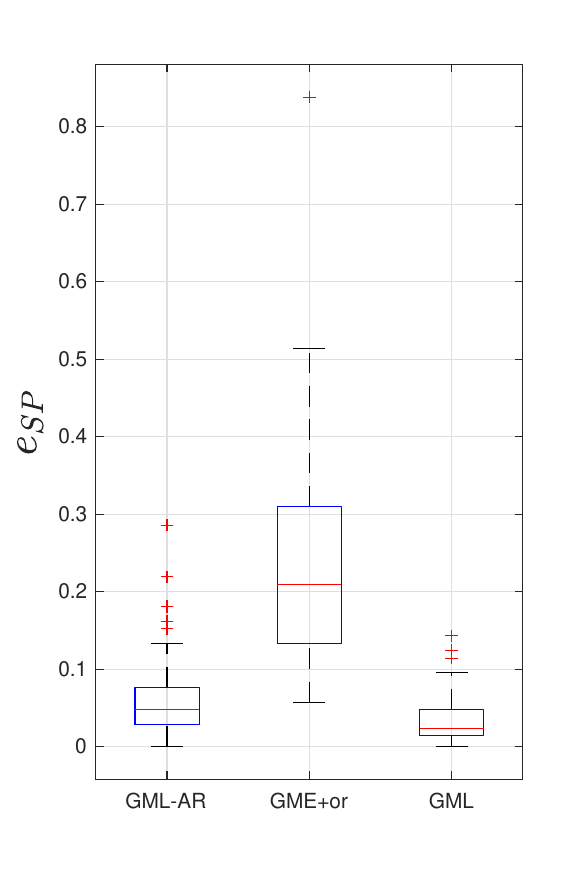}
\caption{Second Monte Carlo experiment with $N=500$. Boxplots of the relative error (left) and the fraction of misspecified edges (right).}\label{figN1000}
\end{figure}

We aim to test the performance of the proposed estimator in respect to the methods available in the literature through some Monte Carlo experiments. Each experiment is composed by $M=100$ trials and in each trial we randomly generate an ARMA graphical model with power spectral density $\Phi_T=\pb_T\Qb_T^{-1}$ with $m=15$ nodes and order $n=2$.   
More precisely, the matrix coefficients of $\Qb_T$ are generated in such a way that the fraction of nonnull entries in $\Qb_T$ is equal to 0.17. The coefficients of $\pb_T$ are generated in such a way that $\pb_T$ has a stable zero whose modulus is greater than 0.98. In this way, the resulting model will not be well approximated by a low-order AR model. Finally, we generate a data sequence $y^N$ of length $N$ from such model. Then, we estimate the model from $y^N$ using the following estimators:
\begin{itemize}
\item \textbf{ME} which denotes the maximum entropy estimator for ARMA models and whose parameters are given by solving (\ref{dual_pb_r}) with $E=V\times V$; this method represents the multivariate extension to the maximum entropy method proposed in \cite{RKL-16multidimensional}; here $\lambda=1$;
\item \textbf{GML-AR} which denotes the generalized maximum likelihood (GML) method for estimating AR graphical models \cite{ZORZI2019108516}; here $\varepsilon=10^{-4}$ and $\epsilon=10^{-8}$;

\item \textbf{GME+or} which denotes the generalized maximum entropy method to estimate ARMA graphical models, proposed in \cite{YOU2022110319}, and equipped with an ``oracle''; more precisely, the oracle knows the true trigonometric polynomial $\pb_{T}$; here the grid for the values of the regularization parameter $\gamma$ has been chosen as in \cite{YOU2022110319}, that is $\{0, 0.02, 0.04, 0.06, 0.08, 0.1, 0.3, 0.5, 1\}$; the best $\gamma$ is selected according to the BIC criterium; finally, the threshold for the stopping condition is $\epsilon=10^{-8}$;
\item \textbf{GML} which denotes the estimator proposed in Section~\ref{sec:Bayes} with $\lambda=1$, $\varepsilon=10^{-4}$ and $\epsilon=10^{-8}$. 
\end{itemize}
To assess the performance of the estimators we compute the fraction of misspecified edges with respect to the true model:
\al{e_{SP}=\frac{\|\Eb-\hat \Eb\|_0}{m^2}\nn}
where $\Eb=\mathrm{supp}(\Phi_T^{-1})$, $\hat \Eb=\mathrm{supp}(\hat \Phi^{-1})$, $\hat \Phi$ denotes the estimator of $\Phi_T$ and $\|\cdot\|_0$ denotes the $\ell_0$ matrix norm. Moreover, we compute the relative error of $\hat \Phi^{-1}$  with respect to the true one $\Phi_T^{-1}$:
\al{err=\frac{\int_{\Tbb} \|\hat \Phi^{-1}-\Phi_T^{-1}\|_F}{\int_{\Tbb} \|\Phi_T^{-1}\|_F}\nn}
where $\|\cdot\|_F$ denotes the Frobenius norm. In the first Monte Carlo experiment we have considered $N=500$. Figure \ref{figN500}(left) shows the boxplots of $err$ for the different estimators. \textbf{GML-AR} exhibits the worst performance because it is not able to approximate the ARMA model through a low order  AR model. \textbf{ME} performs worse than \textbf{GME+or} and \textbf{GML} because it estimates a full graphical model, i.e. this method is not able to control the model complexity in terms of number of edges. Finally, \textbf{GME+or} exhibits the best performance thanks to its oracle. Figure \ref{figN500}(right) shows the boxplots of $e_{SP}$ for \textbf{GML-AR}, \textbf{GME+or} and \textbf{GML}. Notice that, we did not plot the boxplot for \textbf{ME} because it always estimates a full graph, i.e. the fraction of misspecified edges is close to 1.   \textbf{GME+or} exhibits the worst performance because the regularization term used to induce sparsity is 
\al{\gamma \sum_{j\geq h} q_{jh}(\Qb)\nn }  
which penalizes the elelemts of $\Qb$ in the same way. Thus, it introduces an adverse bias in the case $q_{jh}(\Qb)$s take very different values. In addition, the fact that the grid is fixed (otherwise the user should change it for any trial)  could degradate the performance of the estimator.     
Finally, \textbf{GML} performs slightly better than \textbf{GML-AR}. Their performance is similar because they use the same regularization term to induce sparsity. 
In the second Monte Carlo experiment we have considered $N=1000$. Figure \ref{figN1000} shows the boxplots of $err$ and $e_{SP}$ for the different estimators. In regard to the relative error, the performance of \textbf{ME} is improved in respect to the case $N=500$ because increasing the data length the overfitting phenomenon is less evident. However, \textbf{ME} is still worse than \textbf{GML+or} and \textbf{GML}.  Finally,  \textbf{GML} is able to achieve a performance similar to the one of \textbf{GME+or}. In regard to the fraction of misspecified edges, the scenario is similar to the case $N=500$, but the performance for all the estimators now is slightly improved.  To conclude, these two Monte Carlo experiments show the superiority of the proposed estimator in the case the data is generated by an ARMA model. 

\section{Conclusions}\label{sec:concl}
We have introduced a maximum entropy covariance and cepstral extension problem for graphical models. Under mild assumptions, we showed that the problem admits an approximate solution which represents an ARMA graphical model whose topology is determined by the selected entries of the covariance lags considered in the extension problem. Then, we showed how its dual problem is linked with the maximum likelihood estimator.  
From such connection we have proposed a Bayesian model and an approximate maximum a posteriori estimator to learn ARMA graphical models in the case the topology is not known. Finally, we have conducted some numerical experiments to show the effectiveness of the proposed method. It is worth noting the concepts introduced in this paper can be extended to other types of ARMA graphical models, e.g. the latent-variable graphical models, \cite{TAC19}, for which only the AR case has been addressed.


\begin{thebibliography}{10}
\providecommand{\url}[1]{#1}
\csname url@samestyle\endcsname
\providecommand{\newblock}{\relax}
\providecommand{\bibinfo}[2]{#2}
\providecommand{\BIBentrySTDinterwordspacing}{\spaceskip=0pt\relax}
\providecommand{\BIBentryALTinterwordstretchfactor}{4}
\providecommand{\BIBentryALTinterwordspacing}{\spaceskip=\fontdimen2\font plus
\BIBentryALTinterwordstretchfactor\fontdimen3\font minus
  \fontdimen4\font\relax}
\providecommand{\BIBforeignlanguage}[2]{{%
\expandafter\ifx\csname l@#1\endcsname\relax
\typeout{** WARNING: IEEEtran.bst: No hyphenation pattern has been}%
\typeout{** loaded for the language `#1'. Using the pattern for}%
\typeout{** the default language instead.}%
\else
\language=\csname l@#1\endcsname
\fi
#2}}
\providecommand{\BIBdecl}{\relax}
\BIBdecl

\bibitem{e20010076}
S.~Maanan, B.~Dumitrescu, and C.~Giurc{\"a}eneanu, ``Maximum entropy
  expectation-maximization algorithm for fitting latent-variable graphical
  models to multivariate time series,'' \emph{Entropy}, vol.~20, no.~1, 2018.

\bibitem{you2023sparse}
J.~You and C.~Yu, ``Sparse plus low-rank identification for dynamical
  latent-variable graphical {AR} models,'' \emph{arXiv preprint
  arXiv:2307.11320}, 2023.

\bibitem{LATENTG}
M.~Zorzi and R.~Sepulchre, ``{AR} identification of latent-variable graphical
  models,'' \emph{IEEE Trans. on Automatic Control}, vol.~61, no.~9, pp.
  2327--2340, 2016.

\bibitem{KRON}
M.~Zorzi, ``Autoregressive identification of {K}ronecker graphical models,''
  \emph{Automatica}, vol. 119, p. 109053, 2020.

\bibitem{MATERASSI2}
D.~{Materassi} and M.~V. {Salapaka}, ``Signal selection for estimation and
  identification in networks of dynamic systems: A graphical model approach,''
  \emph{IEEE Transactions on Automatic Control}, vol.~65, no.~10, pp.
  4138--4153, 2020.

\bibitem{9600870}
M.~S. Veedu, H.~Doddi, and M.~V. Salapaka, ``Topology learning of linear
  dynamical systems with latent nodes using matrix decomposition,'' \emph{IEEE
  Transactions on Automatic Control}, vol.~67, no.~11, pp. 5746--5761, 2022.

\bibitem{SHAIKHVEEDU2023111182}
M.~{Shaikh Veedu} and M.~V. Salapaka, ``Topology identification under spatially
  correlated noise,'' \emph{Automatica}, vol. 156, p. 111182, 2023.

\bibitem{doddi2020estimating}
H.~Doddi, D.~Deka, S.~Talukdar, and M.~Salapaka, ``Estimating linear dynamical
  networks of cyclostationary processes,'' \emph{arXiv preprint
  arXiv:2009.12667}, 2020.

\bibitem{HOF1}
K.~R. Ramaswamy, G.~Bottegal, and P.~M. {Van den Hof}, ``Learning linear
  modules in a dynamic network using regularized kernel-based methods,''
  \emph{Automatica}, vol. 129, p. 109591, 2021.

\bibitem{HOF2}
K.~R. Ramaswamy and P.~M. J. V.~d. Hof, ``A local direct method for module
  identification in dynamic networks with correlated noise,'' \emph{IEEE
  Transactions on Automatic Control}, vol.~66, no.~11, pp. 5237--5252, 2021.

\bibitem{FPR-08}
A.~Ferrante, M.~Pavon, and F.~Ramponi, ``Hellinger versus {K}ullback--{L}eibler
  multivariable spectrum approximation,'' \emph{IEEE Transactions on Automatic
  Control}, vol.~53, no.~4, pp. 954--967, 2008.

\bibitem{RFP-09}
F.~Ramponi, A.~Ferrante, and M.~Pavon, ``A globally convergent matricial
  algorithm for multivariate spectral estimation,'' \emph{IEEE Transactions on
  Automatic Control}, vol.~54, no.~10, pp. 2376--2388, 2009.

\bibitem{FRT-11}
A.~Ferrante, F.~Ramponi, and F.~Ticozzi, ``On the convergence of an efficient
  algorithm for {K}ullback--{L}eibler approximation of spectral densities,''
  \emph{IEEE Transactions on Automatic Control}, vol.~56, no.~3, pp. 506--515,
  2011.

\bibitem{zhu2018wellposed}
B.~Zhu, ``On the well-posedness of a parametric spectral estimation problem and
  its numerical solution,'' \emph{IEEE Transactions on Automatic Control},
  vol.~65, no.~3, pp. 1089--1099, 2020.

\bibitem{zhu2020m}
B.~Zhu, A.~Ferrante, J.~Karlsson, and M.~Zorzi, ``M$^{2}$-spectral estimation:
  A flexible approach ensuring rational solutions,'' \emph{SIAM Journal on
  Control and Optimization}, vol.~59, no.~4, pp. 2977--2996, 2021.

\bibitem{ringh2018multidimensional}
A.~Ringh, J.~Karlsson, and A.~Lindquist, ``Multidimensional rational covariance
  extension with approximate covariance matching,'' \emph{SIAM Journal on
  Control and Optimization}, vol.~56, no.~2, pp. 913--944, 2018.

\bibitem{KLR-16multidimensional}
J.~Karlsson, A.~Lindquist, and A.~Ringh, ``The multidimensional moment problem
  with complexity constraint,'' \emph{Integral Equations and Operator Theory},
  vol.~84, no.~3, pp. 395--418, 2016.

\bibitem{RKL-16multidimensional}
A.~Ringh, J.~Karlsson, and A.~Lindquist, ``Multidimensional rational covariance
  extension with applications to spectral estimation and image compression,''
  \emph{SIAM Journal on Control and Optimization}, vol.~54, no.~4, pp.
  1950--1982, 2016.

\bibitem{enqvist2004aconvex}
P.~Enqvist, ``A convex optimization approach to {ARMA}$(n,m)$ model design from
  covariance and cepstral data,'' \emph{SIAM Journal on Control and
  Optimization}, vol.~43, no.~3, pp. 1011--1036, 2004.

\bibitem{SONGSIRI_TOP_SEL_2010}
J.~Songsiri and L.~Vandenberghe, ``Topology selection in graphical models of
  autoregressive processes,'' \emph{J. Mach. Learning Res.}, vol.~11, pp.
  2671--2705, 2010.

\bibitem{MAANAN2017122}
S.~Maanan, B.~Dumitrescu, and C.~Giurc{\"a}eneanu, ``Conditional independence
  graphs for multivariate autoregressive models by convex optimization:
  Efficient algorithms,'' \emph{Signal Processing}, vol. 133, pp. 122--134,
  2017.

\bibitem{ALPAGO_SL_REC}
D.~Alpago, M.~Zorzi, and A.~Ferrante, ``A scalable strategy for the
  identification of latent-variable graphical models,'' \emph{IEEE Transactions
  on Automatic Control}, vol.~67, no.~7, pp. 3349--3362, 2022.

\bibitem{6365751}
E.~Avventi, A.~G. Lindquist, and B.~Wahlberg, ``{ARMA} identification of
  graphical models,'' \emph{IEEE Transactions on Automatic Control}, vol.~58,
  no.~5, pp. 1167--1178, 2013.

\bibitem{LINKPRED}
D.~Alpago, M.~Zorzi, and A.~Ferrante, ``Data-driven link prediction over
  graphical models,'' \emph{IEEE Transactions on Automatic Control}, vol.~68,
  no.~4, pp. 2215--2228, 2023.

\bibitem{YOU2022110319}
J.~You, C.~Yu, J.~Sun, and J.~Chen, ``Generalized maximum entropy based
  identification of graphical {ARMA} models,'' \emph{Automatica}, vol. 141, p.
  110319, 2022.

\bibitem{WIPF_2010}
D.~Wipf and S.~Nagarajan, ``Iterative reweighted $\ell_1$ and $\ell_2$ methods
  for finding sparse solutions,'' \emph{J. Sel. Topics Signal Processing},
  vol.~4, no.~2, pp. 317--329, 2010.

\bibitem{ZORZI2019108516}
M.~Zorzi, ``Empirical {B}ayesian learning in {AR} graphical models,''
  \emph{Automatica}, vol. 109, p. 108516, 2019.

\bibitem{brillinger1996}
D.~R. Brillinger, ``Remarks concerning graphical models for time series and
  point processes,'' \emph{Brazilian Review of Econometrics}, vol.~16, no.~1,
  pp. 1--23, 1996.

\bibitem{lindquist2013multivariate}
A.~Lindquist, C.~Masiero, and G.~Picci, ``On the multivariate circulant
  rational covariance extension problem,'' in \emph{CDC}.\hskip 1em plus 0.5em
  minus 0.4em\relax IEEE, 2013, pp. 7155--7161.

\bibitem{DEMPSTER_1972}
A.~Dempster, ``Covariance selection,'' \emph{Biometrics}, vol.~28, no.~1, pp.
  157--175, 1972.

\bibitem{Zhu-Zorzi2023cepstral}
B.~Zhu and M.~Zorzi, ``A well-posed multidimensional rational covariance and
  generalized cepstral extension problem,'' \emph{SIAM Journal on Control and
  Optimization}, vol.~61, no.~3, pp. 1532--1556, 2023.

\bibitem{rudin1987real}
W.~Rudin, \emph{Real and Complex Analysis}, 3rd~ed.\hskip 1em plus 0.5em minus
  0.4em\relax McGraw-Hill Book Company, 1987.

\bibitem{ringh2015multidimensional}
A.~Ringh, J.~Karlsson, and A.~Lindquist, ``The multidimensional circulant
  rational covariance extension problem: Solutions and applications in image
  compression,'' in \emph{54th Annual Conference on Decision and Control
  (CDC)}.\hskip 1em plus 0.5em minus 0.4em\relax IEEE, 2015, pp. 5320--5327.

\bibitem{FALCONI2023110672}
L.~Falconi, A.~Ferrante, and M.~Zorzi, ``Mean-square consistency of the
  f-truncated {M}$^2$-periodogram,'' \emph{Automatica}, vol. 147, p. 110672,
  2023.

\bibitem{whittle1953estimation}
P.~Whittle, ``Estimation and information in stationary time series,''
  \emph{Arkiv f{\"o}r matematik}, vol.~2, no.~5, pp. 423--434, 1953.

\bibitem{RASMUSSEN_WILLIAMNS_2006}
C.~Rasmussen and C.~Williams, \emph{{Gaussian Processes for Machine
  Learning}}.\hskip 1em plus 0.5em minus 0.4em\relax The MIT Press, 2006.

\bibitem{zhou1997approximate}
Z.~Zhou, R.~Leahy, and J.~Qi, ``Approximate maximum likelihood hyperparameter
  estimation for gibbs priors,'' \emph{IEEE transactions on image processing},
  vol.~6, no.~6, pp. 844--861, 1997.

\bibitem{TAC19}
V.~Ciccone, A.~Ferrante, and M.~Zorzi, ``Learning latent variable dynamic
  graphical models by confidence sets selection,'' \emph{IEEE Trans. Autom.
  Control}, vol. accepted, 2020.

\end{thebibliography}
\end{document}